\documentclass[a4,12pt,multicol,makeidx]{article}
\def\origin{
  \clearpage
\vskip-\baselineskip\vskip-\topskip%
  \vbox to 0pt{\vskip-1in%
    \hbox to 0pt{\hskip-1in%
      \hbox to 0pt{\vrule width 1cm height .4pt depth 0mm\hss}%
      \vbox to 0pt{\hrule width .4pt height 0pt depth 1cm\vss}%
    \hss}%
  \vss}
  \vskip-\baselineskip
  \vbox to 0pt{\vskip-1in\vskip3cm%
    \hbox to 0pt{\hskip-1in\hskip3cm%
      \hbox to 0pt{\hss\vrule width 2cm height .4pt depth 0mm\hss}%
      \vbox to 0pt{\vss\hrule width .4pt height 1cm depth 1cm\vss}%
    \hss}%
  \vss}%
\vskip5mm\hskip10mm (3cm,3cm)
}%
   
\def\a{\alpha}   
\def\ue{u_{\varepsilon}}
\def\ou{\overline{u}} \def\oI{\overline{I}} \def\ox{\overline{x}} \def\oy{\overline{y}}
\def\ov{\overline{v}} \def\uv{\underline{v}} \def\ow{\overline{w}} \def\uw{\underline{w}}
\def\hx{\hat{x}} \def\hy{\hat{y}} \def\oy{\overline{y}}
\def\ov{\overline{v}}   \def\ty{\tilde{y}}
 \def\l{\lambda}   \def\p{\partial}   
    
\def\leq{\underline{<}} 
\def\la{\langle} \def\ra{\rangle} \def\hx{\hat{x}} \def\hy{\hat{y}}
 \def\I{{\bf1}_{|z|\leq 1}}

\newenvironment{theorem}{%
\par \bigskip \it}{%
\bigskip \par}

\newenvironment{definition}{%
\par \bigskip \it}{%
\bigskip \par}

\makeindex
\title{Quasi-periodic and almost periodic homogenizations of integro-differential equations with L{\'e}vy operators. 
}
\author{Mariko Arisawa\\ DAMTP
\\University of Cambridge\\
E-mail: M.Arisawa@damtp.cam.ac.uk
}
\date{}
\begin{document}
\maketitle
\bigskip
\section{Introduction.} 

$\qquad$ The quasi-periodic homogenization  and the almost periodic homogenization of a class of  integro-differential equations with  L{\'e}vy operators are studied 
in this paper. First, the quasi-periodic homogenization  is the following. For $\varepsilon=$$(\varepsilon_1,\varepsilon_2)\in {\bf R^+}\times {\bf R^+}$, consider
$$
	\ue(x)+ \sup_{\a\in \cal{A}} \{
 	\left\langle -b(x,\a), \nabla\ue \right\rangle \} - a(\frac{x}{\varepsilon_1})
	\int_{{\bf R^N}} [\ue(x+z)-\ue(x)
$$
\begin{equation}\label{qpex}
	- \I \left\langle z, \nabla\ue(x) \right\rangle] \frac{1}{|z|^{N+\a}}dz
	-g_1(\frac{x}{\varepsilon_1})-g_2(\frac{x}{\varepsilon_2})=0 \qquad x\in {\Omega}, 
\end{equation}
\begin{equation}\label{bc}
	\ue(x)=h(x) \qquad x\in {\Omega^c}, 
\end{equation}
where $\Omega$ is an open domain in ${\bf R^N}$, $\cal A$ is a compact subset of a metric space, $a(\cdot)$ a bounded continuous function, 
 $b(x,\a)$  is a bounded function from ${\bf R^N}\times \cal A$ to ${\bf R^N}$ such that there 
exists a constant $L>0$ 
\begin{equation}\label{b}
	|b(x,\a)-b(y,\a)|\leq L|x-y| \qquad \forall x,y\in {\bf R^N},\quad \a\in 
	{\mathcal A}, 
\end{equation}
 $a(y)$, $g_i(y)$ $(i=1,2)$ are real valued periodic functions in ${\bf T^N}$ ($N$ dimensional torus with periods 1), two parameters $\varepsilon_1$, $\varepsilon_2$ 
 satisfy
\begin{equation}\label{irratio}
	\frac{\varepsilon_2}{\varepsilon_1}=\gamma\in {\bf R}\backslash{\bf Q}, 
\end{equation} 
the nonlocal (integral) term is the L{\'e}vy operator with the $\alpha$-stable symmetric measure
$$
	\frac{1}{|z|^{N+\a}}dz \qquad \hbox{with}\quad \a\in (0,2), 
$$
 and $h(x)$ is a bounded  continuous function defined in $\Omega^c \subset {\bf R^N}$. 
The existence and the uniqueness of the solution  $\ue$ is known in the framework of the viscosity solution. 
 We are interested in the asymptotic limit of $\ue$ as $\varepsilon\to 0$ while satisfying the relationship (\ref{irratio}). 
 We shall show in below in more generality the unique existence of the limit 
$\ou=\lim_{\varepsilon_1,\varepsilon_2\to 0} \ue$ and its  characterization by  an effective integro-differential 
equation.  \\

$\qquad$Next, the following is an example of the almost periodic homogenization. Let $\varepsilon>0$, and consider 
$$
	\ue+ \sup_{\a\in \cal{A}} \{
 	\left\langle -b(x,\a), \nabla\ue \right\rangle \} - a(\frac{x}{\varepsilon}) \int_{{\bf R^N}} [\ue(x+z)-\ue(x)
	\qquad\qquad\qquad\qquad\qquad
$$
\begin{equation}\label{almostex}
	- \I
	\left\langle z, \nabla\ue(x) \right\rangle] \frac{1}{|z|^{N+\a}}dz
	-g(\frac{x}{\varepsilon})=0 \qquad x\in {\bf R^N}
\end{equation}
with (\ref{bc}),  
where $a(y)$, $g(y)$ are real valued functions defined in ${\bf R^N}$, uniformly almost periodic in the sense of Bohr \cite{bohr}.  \\

(Uniformly almost periodiic function) A real valued function  $f(y)$ defined in ${\bf R^N}$ is uniformly almost periodic in the sense of Bohr if and only if the set 
 of functions 
$$
	\{f(y+z)| \quad z\in {\bf R^N} \}
$$
is relatively compact in the space of the bounded functions in ${\bf R^N}$ with the norm $\|f \|_{\infty}$$=\sup_{x\in {\bf R^N}}|f(x)|$. \\

Remark that a quasi-periodic function (for example $g_1(\frac{x}{\varepsilon_1})+g_2(\frac{x}{\varepsilon_2})$ in (\ref{qpex})) is a uniformly almost periodic function. 
We refer the readers to  Besicovich \cite{besicovitch} for the rich informations on the uniformly almost periodic function, some of which 
we utilize in below. 
As before, we are interested in the asymptotic limit of the solution $\ue$ of (\ref{almostex}) as $\varepsilon\to 0$, and in characterizing the limit by finding an effective integro-differential equation for it.\\

$\quad$
The present work is a straight forward generalization of the periodic homogenization for the integro-differential equation with the L{\'e}vy operator in Arisawa \cite{ar8} 
to the  quasi-periodic and  the almost periodic homogenizations. In the case of the partial differential equation (PDE in short), such generalizations were done  in 
 Arisawa \cite{ar1} (the quasi-periodic homogenization for first-order PDEs with non-convex Hamiltonians), \cite{ar2} (the almost periodic homogenization for  second-order elliptic PDEs), and in Ishii \cite{i} (the almost periodic homogenization for  first-order PDEs with the convex Hamiltonian), and then more generally treated in the stationally ergodic setting by Caffarelli, Souganidis and Wang in \cite{caff} (the stochastic homogenization for second-order uniformly elliptic PDEs).  It is known that for non-convex first-order Hamilton-Jacobi equations,  the almost-periodic homogenization is not well-posed in general. It is also known in Lions and Souganidis \cite{lp} that the stochastic homogenization for the first-order Hamilton-Jacobi equation is not necessarily well-posed.  
In the integro-differential problems (\ref{qpex}) and (\ref{almostex}),  the $\a$-stable L{\'e}vy operator is  the fractional power of Laplacian: $\Delta^{\frac{\a}{2}}$.    According to wheather $\alpha \in (0,1]$ or 
 $\a\in (1,2)$, the operator can be considered to be close to the  first-order operator or to the second-order elliptic operator. 
 Therefore, the quasi-periodic and the almost periodic homogenizations are natural  to be studied for the integro-differential equations. \\

Now, let us explain the outline of this paper.  We generalize the quasi-periodic problem to the folliwing. 
\begin{equation}\label{quasi}
	\ue+ \sup_{\a\in \cal{A}} \{
 	\left\langle -b(x,\a), \nabla\ue \right\rangle \} - a(\frac{x}{\varepsilon_1})\int_{{\bf R^N}} [\ue(x+z)-\ue(x)
\end{equation}
$$
	- \I \left\langle z, \nabla\ue(x) \right\rangle] \frac{1}{|z|^{N+\a}}dz
	-g_M(\frac{x}{\varepsilon_1},...,\frac{x}{\varepsilon_M})=0 \qquad x\in {\bf R^N}, 
$$
where  such that 
\begin{equation}\label{a}
	a(y_1)\quad\hbox{ is periodic  in}\quad  {\bf T^N};\quad a(y_1)\geq a_0 \qquad \forall y\in {\bf T^N}, 
\end{equation} 
where $a_0>0$ is a constant, 
and $g_M(y_1,...,y_M)$ $(M\in {\bf N})$ is a real valued periodic functions in $(y_1,...,y_M)\in $${\bf T^{MN}}$,  $\varepsilon_i>0$  $(1\leq i\leq M)$ satisfy the following 
 non-resonance condition. \\

(Non-resonance condition) A countable set of real numbers $E=$$\{\varepsilon_i\}$ $(i\in {\bf N})$ is said to satisfy the non-resonance condition if  for any $k\in {\bf N}$ 
 and for $E_{k}=\{ \varepsilon_1,\varepsilon_2,...,\varepsilon_k \}$ the only rational numbers $a_1$, $a_2$,...,$a_k$ to satisfy 
\begin{equation}\label{nonreso}
	\sum_{i=1}^{k} a_i \varepsilon_i =0
\end{equation}
are $a_i=0$ ($1\leq \forall i\leq k$). \\

The above condition is taken from \cite{{besicovitch}}, where the finite version was used in  Arisawa and Lions \cite{al}. 
We assume also that there exists a constant $\theta_0\in (0,1]$ such that 
\begin{equation}\label{holder}
	|a(y)-a(y')|  \leq C|y-y'|^{\theta_0}   \quad \forall y,y'\in {\bf T^N}, \quad 1\leq \forall i\leq M, 
\end{equation}
$$
	|g_M(y_1,...y_{i-1},\oy_i,y_{i+1}...,y_M)-g_M(y_1,...y_{i-1},\oy'_i,y_{i+1}...,y_M)|\leq C|\oy_i-\oy'_i|^{\theta_0}  
$$
\begin{equation}\label{gholder}
	\qquad\qquad\qquad\qquad\qquad\qquad\qquad\qquad
	\quad \forall \oy_i,\oy'_i \in {\bf T^N}, \quad 1\leq \forall i\leq M, 
\end{equation}
where $C>0$ is a constant which depends only on $\theta_0$. \\

 Our method is based on the relationship between the   formal asymptotic expansion and the ergodic problem. 
  The   formal asymptotic expansion was introduced by Bensoussan, Lions and Papanicolaou in \cite{blp},  and developped rigorously by Lions, Papanicolaou and Varadhan \cite{lpv}, Evans \cite{ev1}, \cite{ev2}, and others. In \S 2,  we utilize the formal asymptotic expansion method to 
 obtain the ergodic cell problem. The key ingredient to solve the ergodic cell problem is 
 the strong maximum principle for the L{\'e}vy operator. In  \S 3, we prove the strong maximum principle for some general class of L{\'e}vy operators
$$
	 \int_{{\bf R^N}} [u(x+z)-u(x)- \I
	\left\langle z, \nabla u(x) \right\rangle] dq(z), 
$$
 which includes  the $\a$-stable symmetric operators as special cases. In \S 4, by using the result in \S 3  the quasi-periodic ergodic cell problems are solved. 
 In \S 5, the almost priodic ergodic cell problems are solved.  In \S 6, we give our main results on the  quasi-periodic and the almost periodic homogenizations. \\

For an upper semi-continuous (USC in short) function $u$ and a lower semi-continuous (LSC in short) function $v$ in ${\bf R^N}$, $J^{2,+}_{\Omega}u(x)$ and $J^{2,-}_{\Omega}v(x)$ represent respectively 
 the set of second-order subdifferentials and the set of superdifferentials of $u$ and $v$ at $x \in \Omega$. That is, for $u\in USC({\bf R^N})$, $(p,Q)\in J^{2,+}_{\Omega}u(x)$ 
means that $(p,Q)\in {\bf R^N}\times{\bf S^N}$, and for any $\delta>0$ there exists $\nu>0$ such that
\begin{equation}\label{pX}
	u(x+z)\leq u(x)+\la p,z \ra+ \frac{1}{2} \la Qz,z \ra +\delta |z|^2 \quad \forall |z|\leq \nu. 
\end{equation}
For $v\in LSC({\bf R^N})$,  $(p,Q)\in J^{2,-}_{\Omega}v(x)$ means that $(p,Q)\in {\bf R^N}\times{\bf S^N}$, and for any $\delta>0$ there exists $\nu>0$ such that 
\begin{equation}\label{pX2}
	v(x+z)\geq v(x)+\la p,z \ra+ \frac{1}{2} \la Qz,z \ra -\delta |z|^2  \quad \forall |z|\leq \nu. 
\end{equation}
We use the notation $I[u](x)$$=$$\int_{{\bf R^N}}[u(x+z)-u(x)-\I \la\nabla u(x),z \ra] dq(z)$, 
$$
	I_{\nu,\delta}^{1,+}[u,p,X](x) = \int_{|z|\leq \nu}\frac{1}{2}\la(X+2\delta I)z,z \ra dq(z),
$$
(resp.
$$
	I_{\nu,\delta}^{1,-}[u,p,X](x) = \int_{|z|\leq \nu}\frac{1}{2}\la(X-2\delta I)z,z \ra dq(z),)
$$
$$
	I_{\nu,\delta}^{2}[u,p,X](x) = \int_{|z|> \nu} [u(x+z)-u(x)- {\bf 1}_{|z|\leq 1}\la p,z \ra] dq(z).
$$
 Consider 
\begin{equation}\label{A}
	A(x,u(x),\nabla u(x),\nabla^2 u(x),I[{u}](x))=0
	\quad x\in {\Omega}, 
\end{equation}
where $A(x,u,p,Q,I)$$\in C(\Omega \times {\bf R}\times{\bf R^N} \times {\bf S^N} \times {\bf R} )$. 

\begin{definition}{\bf Definition 1.1.}  Let $u\in USC({\bf R^N})$ (resp. $v\in LSC({\bf R^N})$). We say that $u$ (resp. $v$) is a viscosity subsolution (resp. supersolution) of (\ref{A}), if for any $\hx\in \Omega$, any $(p,X)\in J_{{\bf R^N}}^{2,+}u(\hx)$ (resp. $J_{{\bf R^N}}^{2,-}v(\hx)$), and 
for any pair of numbers $(\varepsilon,\delta)$ satisfying  (\ref{pX})  (resp. (\ref{pX2})), 
the following holds 
$$
	A(\hat{x},u(\hat{x}),p,X, I_{\nu,\delta}^{1,+}[u,p,X](\hx)
	+ I_{\nu,\delta}^{2}[u,p,X](\hx)) \leq 0.
$$
(resp.
$$
	A(\hat{x},v(\hat{x}),p,X, I_{\nu,\delta}^{1,-}[v,p,X](\hx)
	+ I_{\nu,\delta}^{2}[v,p,X](\hx)) \geq 0.
$$
)
 If $u$ is a viscosity subsolution and a viscosity supersolution at the same time, it is called a viscosity solution.
\end{definition}

 The above Definition 1.1 is equivalent to the following (see Arisawa \cite{ar7}). \\

\begin{definition}{\bf Definition 1.2.}  Let $u\in USC({\bf R^N})$ (resp. $v\in LSC({\bf R^N})$). We say that $u$ (resp. $v$) is a viscosity subsolution (resp. supersolution) of (\ref{A}),  if for any $\hx\in \Omega$, any $\phi\in C^2({\bf R^N})$ such that $u(\hx)=\phi(\hx)$ and $ u-\phi $ takes a global maximum (resp. minimum) at $\hx$, 
\begin{equation}\label{defb1}
	A(\hat{x},u(\hat{x}),\nabla \phi(\hat{x}),\nabla^2 \phi(\hat{x}), I[\phi](\hx)) \leq 0.
\end{equation}
(resp.
\begin{equation}\label{defb2}	
	A(\hat{x},v(\hat{x}),\nabla \phi(\hat{x}),\nabla^2 \phi(\hat{x}), I[\phi](\hx)) \geq 0.)
\end{equation}
If $u$ is a viscosity subsolution and a viscosity supersolution at the same time, it is called a viscosity solution.
\end{definition}

The existence and the uniqueness of the solution $\ue$ of (\ref{almostex})-(\ref{bc}) and (\ref{quasi})-(\ref{bc})  are established in the framework of the viscosity solution. We refer the readers to Arisawa \cite{ar3}, \cite{ar5}, Barles, Buckdahn and Pardoux \cite{bbp}, Barles and Imbert \cite{bi}, etc...\\

\section{Formal asymptotic expansions.}

$\qquad$We treat  the quasi-periodic homogenization (\ref{quasi}). The almost periodic homogenization can be treated similarly, which we mention in \S 5 below. Put $\varepsilon$$=(\varepsilon_1,...,\varepsilon_M)$, and $\gamma_i=\frac{\varepsilon_i}{\varepsilon_1}$ ($1\leq i\leq M$). We devide the situation into  three cases. \\

\begin{itemize}
\item I. $\a\in (0,1)$ and $b(x,\a)\not \equiv 0$.
\item II. $\a=1$ and $b(x,\a)\not \equiv 0$.
\item III. $\a\in (1,2)$, or $\a\in (0,1]$ and $b(x,\a)\equiv 0$.
\end{itemize}

The formal asymptotic expansions are: for the case of I and II 
\begin{equation}\label{formal1} 
	\ue(x)=\ou(x) + \varepsilon_1 v(\frac{x}{\varepsilon_1}), 
\end{equation}
and for the case of III
\begin{equation}\label{formal2} 
	\ue(x)=\ou(x) + \varepsilon_1^{\a}v(\frac{x}{\varepsilon_1}). 
\end{equation}

We introduce in (\ref{quasi}) the formal derivatives of the above expansions for each case. \\

{\bf Case I.} If  $\a\in (0,1)$ and $b(x,\a)\not \equiv 0$,  by introducing the derivatives of (\ref{formal1}) into (\ref{quasi}),  ignoring the $o(1)$ terms, and rewriting 
$y=\frac{x}{\varepsilon_1}$, $p=\nabla \ou(x)$, we get the following relationship. 
$$
	\ou(x)+ \sup_{\a\in \cal{A}} \{
 	\left\langle -b(x,\a), p+ \nabla_y v(y) \right\rangle \} 
	- a(y) I [\ou](x)- g_M(\gamma_1^{-1}y,...,\gamma_M^{-1}y)=0. 
$$
As in \cite{lpv} and other works, for each fixed $(x,p,I)\in$ 
$\Omega \times {\bf R^N}\times {\bf R}$ we intend to get  a unique constant $d_{x,p,I}$
 such that there exists at least a viscosity solution $v(y)$, bounded in  ${\bf R^N}$, 
\begin{equation}\label{cell1}
	d_{x,p,I} +\sup_{\alpha\in \cal{A}}\{\la -b(x,\a),p+\nabla_y v(y) \ra \}- a(y) I -g_M(\gamma_1^{-1}y,...,\gamma_M^{-1}y)=0 \quad y\in {\bf R^N}, 
\end{equation}
which is the ergodic cell problem. 
In fact, if $d_{x,p,I}$ exists for any $(x,p,I)$, then by putting $\overline{I}(x,p,I)=-d_{x,p,I}$ the limit $\ou$ formally satisfies 
\begin{equation}\label{effect}
	\ou+\overline{I}(x,\nabla \ou (x),I[\ou](x))=0 \qquad x\in {\Omega}, 
\end{equation}
which will be verified rigorously  below in \S 6. \\

{\bf Case II.} If $\a=1$ and $b(x,\a)\not \equiv 0$, the introduction of the derivatives of (\ref{formal1}) into (\ref{quasi}) leads
$$
	\ou+ \sup_{\a\in \cal{A}} \{
 	\left\langle -b(x,\a), p+ \nabla_y v(y) \right\rangle \} 
	- a(y)\int_{{\bf R^N}} [v(y+z)-v(y)
$$
$$
	-{\bf 1}_{|z|\leq 1}\left\langle z, \nabla_y v(y) \right\rangle] \frac{1}{|z|^{N+\a}}dz
	-a(y) I[u](x) -g_M(\gamma_1^{-1}y,...,\gamma_M^{-1}y)=0, 
$$
where $y$, $p$ are same as above. The following ergodic cell problem is thus derived. For each fixed $(x,p,I)\in {\Omega \times {\bf R^N} \times {\bf R}}$,  
find a unique number 
 $d_{x,p,I}$ such that the following problem has at least a viscosity solution $v(y)$, bounded in  ${\bf R^N}$,  
\begin{equation}\label{cell2}
	d_{x,p,I} 
	+\sup_{\alpha\in \cal{A}}\{\la -b(x,\a),p+\nabla_y v(y)\ra \}-a(y) \int_{{\bf R^N}} [v(y+z)-v(y)
\end{equation}
$$
	- {\bf 1}_{|z|\leq 1} \left\langle z, \nabla_y v(y) \right\rangle] \frac{1}{|z|^{N+\a}}dz -a(y)I-g_M(\gamma_1^{-1}y,...,\gamma_M^{-1}y)=0 \quad y\in {\bf R^N}. 
$$
As before,  if $d_{x,p,I}$ exists for any $(x,p,I)$, then by putting $\overline{I}(x,p,I)=-d_{x,p,I}$ the limit $\ou$ formally satisfies 
 (\ref{effect}), which will be shown rigorously later. \\

{\bf Case III.} If $\a\in (1,2)$, or $\a\in (0,1]$ and $b(x,\a)\equiv 0$, the introduction of the derivatives of (\ref{formal2}) into (\ref{quasi}) leads 
$$
	\ou+ \sup_{\a\in \cal{A}} \{
 	\left\langle -b(x,\a), \nabla\ou \right\rangle \} 
	- a(y)I[\ou](x)- a(y)
	\int_{{\bf R^N}} [\ou(x+z)-
	\ou(x) \qquad \quad
$$$$
	- \I\left\langle z, \nabla\ou(x) \right\rangle] \frac{1}{|z|^{N+\a}}dz
	- \int_{{\bf R^N}} [v(y+z)-v(y)- 
	\I \left\langle z, \nabla_y v(y) \right\rangle] \frac{1}{|z|^{N+\a}}dz
$$
$$
	\qquad\qquad\qquad\qquad\qquad\qquad\qquad\qquad\qquad
	-g_M(\gamma_1^{-1}y,...,\gamma_M^{-1}y)=0, 
$$
where $y=\frac{x}{\varepsilon_1}$, $p=\nabla \ou(x)$, and $o(1)$ terms are neglected. Then, we are interested in the following ergodic cell problem.  
For each fixed $(x,p,I)$$\in \Omega \times {\bf R^N}\times {\bf R}$, find a unique constant $d_{x,p,I}$
 such that there exists at least a viscosity solution $v(y)$, bounded in  ${\bf R^N}$, 
\begin{equation}\label{cell3}
	d_{x,p,I}- a(y)\int_{{\bf R^N}} [v(y+z)-v(y)
	- \I \left\langle z, \nabla_y v(y) \right\rangle] \frac{1}{|z|^{N+\a}}dz- a(y)I 
\end{equation}
$$
	\qquad\qquad\qquad\qquad\qquad\qquad\qquad\qquad -g_M(\gamma_1^{-1}y,...,\gamma_M^{-1}y)=0 \qquad y\in {\bf R^N}. 
$$
 If $d_{x,p,I}$ exists for any $(x,p,I)$, then by defining $\overline{I}(x,p,I)=-d_{x,p,I}$, 
the limit $\ou$ formally satisfies  (\ref{effect}), 
which will  be rigorously proved  in below.  \\

$\quad$Instead of ({\ref{formal1}}) and ({\ref{formal2}}), the following expansions are also possible. For the cases of I and II 
\begin{equation}\label{formalq1} 
	\ue(x)=\ou(x) + \varepsilon_1 w(\frac{x}{\varepsilon_1},\frac{x}{\varepsilon_2},...,\frac{x}{\varepsilon_M}), 
\end{equation}
and for the case of III
\begin{equation}\label{formalq2} 
	\ue(x)=\ou(x) + \varepsilon_1^{\a} w(\frac{x}{\varepsilon_1},\frac{x}{\varepsilon_2},...,\frac{x}{\varepsilon_M}). 
\end{equation}

Let $B(x,\a)=(\gamma_1^{-1}b(x,\a),...,\gamma_M^{-1}b(x,\a))$, $\Gamma z=(\gamma_{1}^{-1}z,...\gamma_{M}^{-1}z)$. 
By introducing the derivatives of (\ref{formalq1}), (\ref{formalq2}) into (\ref{quasi}), by putting $p=\nabla_x \ou$, $y_i=\frac{x}{\varepsilon_i}$, $I=I[\ou](x)$, we get the ergodic cell problems:  find a unique number $d_{x,p,I}$ such that there exists at least a periodic viscosity solution $w(\oy)$ ($\oy=(y_1,...,y_M)$$\in {\bf T^{MN}}$,  $y_i\in {\bf T^N}$, $1\leq \forall i\leq M$) which satisfies the following. 
 For the case I, 
\begin{equation}\label{qcell1}
	d_{x,p,I} +\sup_{\alpha\in \cal{A}}\{\la -b(x,\a),p\ra - \la B(x,\a),\nabla w(\oy) \ra \}- a(y_1) I -
	g_M(\oy)=0 \quad \oy\in {\bf T^{MN}}. 
\end{equation}
For the case II, 
$$
	d_{x,p,I} +\sup_{\alpha\in \cal{A}}\{\la -b(x,\a),p\ra - \la B(x,\a),\nabla w(\oy) \ra \}- a(y_1) I -a(y_1) \int_{{\bf R^N}} [w(\oy+\Gamma z)
$$
\begin{equation}\label{qcell2}\qquad \qquad
	-w(\oy)- {\bf 1}_{|z|\leq 1} \left\langle \Gamma z, \nabla w(\oy) \right\rangle] dq(z)
	- g_M(\oy)=0 \quad \oy\in {\bf T^{MN}}. 
\end{equation}
For the case III,
 $$
	d_{x,p,I}  -a(y_1) \int_{{\bf R^N}} [w(\oy+\Gamma z)-w(\oy) - {\bf 1}_{|z|\leq 1} \left\langle \Gamma z, \nabla w(\oy) \right\rangle] dq(z)- a(y_1) I
$$
\begin{equation}\label{qcell3}
	\qquad\qquad\qquad\qquad\qquad\qquad\qquad\qquad\qquad
	-g_M(\oy)=0 \quad \oy\in {\bf T^{MN}}. 
\end{equation}

{\bf Remark 2.1.}  The two types of the ergodic cell problems (\ref{cell1}), (\ref{cell2}), (\ref{cell3}) and (\ref{qcell1}), (\ref{qcell2}), (\ref{qcell3}) are respectively connected   by the relationship 
\begin{equation}\label{vw}
	v(y)=w(y,\gamma_2^{-1}y,...,\gamma_M^{-1}y). 
\end{equation}
We use (\ref{qcell1})-(\ref{qcell3}) to complement the informations of (\ref{cell1}),  (\ref{cell2}) and (\ref{cell3}) in below.\\

\section{Strong maximum principle.} 

$\qquad$The strong maximum principle is the key to solve the ergodic cell problem. The present result concerns with a general class of the L{\'e}vy operators 
including the $\a$-stable symmetric operator.  This is an improvement of our previous result in \cite{ar5}.  Consider
\begin{equation}\label{strong}
	H(x,\nabla u,\nabla^2 u)-\int_{{\bf R^N}} [u(x+z)-u(x)-\I \left\langle z, \nabla u(x) \right\rangle ]dq(z)=0\quad x\in {\bf R^N}, 
\end{equation}
where $H$$\in$$C(\Omega\times{\bf R^N}\times{\bf S^N})$, $dq(z)$ is a positive Radon measure  such that 
\begin{equation}\label{Radon} 
	\int_{|z|<1} |z|^2    dq(z)  + \int_{|z|\geq 1} |z|^2    dq(z) < \infty. 
\end{equation}
Assume that 
\begin{equation}\label{positive}
	H(x,0,O)\geq 0 \qquad \forall x\in {\bf R^N}, 
\end{equation}
and that there exists a ball in ${\bf R^N}$, $B=B_r(0)$, centered at the origin with radius $r>0$, for which the following holds 
\begin{equation}\label{trans}
	\int_{B} {\bf 1} dq(z)>0 \qquad \forall x\in {\bf R^N}. 
\end{equation}

{\bf Theorem 3.1.$\quad$} 
\begin{theorem}
Let $u$ be a viscosity subsolution of  (\ref{strong}), and assume that (\ref{positive}) and (\ref{trans}) hold. 
 Assume also that 
there exists a maximum point of $u$, $\hx$ in ${\bf R^N}$. Then, $u$ is constant almost everywhere in ${\bf R^N}$. 
\end{theorem}
${\bf Proof.}$ Let $M=\max_{{\bf R^N}}u(x)$. Put $D=\{x\in {\bf R^N}|\quad u(x)=M\}$, which  is non-empty and closed from the assumption. If $D={\bf R^N}$, the claim is clear. So,  
assume that there exists a point $y_1\in D^c$ such that $dist(y_1,D)=\inf_{x\in D}|y_1-x|\leq \frac{r}{2}$. Take $x_1\in D$ such that $|x_1-y_1|\leq \frac{r}{2}$. 
Since $D^c$ is open, there exists $0<s<\frac{r}{4}$ such that 
\begin{equation}\label{small}
	u(y)< u(x_1)=M \quad \hbox{if}\quad  |y-y_1|< s. 
\end{equation}
Since $x_1$ is a maximum point of $u$, $(0,O)\in J^{2,+}_{{\bf R^N}}u(x_1)$, i.e. for any $\delta>0$ there exists $\nu>0$ such that 
$$
	u(x_{1}+z)\leq u(x_{1}) + \left\langle 0,z  \right\rangle + \frac{1}{2} \left\langle Oz,z  \right\rangle +\delta |z|^2 \quad \forall |z|\leq \nu. 
$$
From the definition of the viscosity subsolution
$$
	H(x_1,0,O) 
	- \int_{|z|\leq \nu}  \frac{1}{2}\left\langle (O+2\delta I)z,z  \right\rangle dq(z) 
	\qquad\qquad\qquad\qquad
$$
$$
	\qquad\qquad\qquad\qquad- \int_{|z|> \nu} [ u(x_1+z)-u(x_1)
	-\I \left\langle 0,z  \right\rangle ]dq(z) \leq 0. 
$$
Put $E=\{z|x_1+z\in B(y_1,s)\}$. Remark that $E\subset B$. Since $u(x_1+z)-u(x_1)\leq 0$ for any $z\in {\bf R^N}$, from (\ref{positive}), (\ref{trans}) and (\ref{small}), the above inequality leads 
$$
	0<- \int_{E\cap \{|z|> \nu\}} [ u(x_1+z)-u(x_1)-\left\langle 0,z  \right\rangle ] dq(z) \qquad
$$
$$
	\leq - \int_{|z|> \nu} [ u(x_1+z)-u(x_1)-\left\langle 0,z  \right\rangle ] dq(z) \leq O(\delta). 
$$
By choosing $\delta>0$ small enough we get  a contradiction, and $D={\bf R^N}$ must hold. 
 \\

{\bf Remark 3.1.} 1. In \cite{ar5}, instead of (\ref{trans}), the following condition was assumed.  
\begin{equation}\label{trans2}
	\int_{D} {\bf 1} dq(z)>0 \qquad \forall x\in {\bf R^N}, \quad \forall D\subset {\bf R^N} \quad \hbox{open}.
\end{equation}
Various generalization is possible beyond Theorem 2.1, which we shall visit in our future work. \\
2. The $\a$-stable symmetric operator satisfies the conditions (\ref{Radon}) and (\ref{trans}) assumed in 
Theorem 3.1.\\

For the later purpose, we are also interested in the 
following "degenerate" L{\'e}vy  operator in ${\bf T^{2N}}$. Let $v(x_1,x_2)$ be a periodic function in ${\bf R^{2N}}$,  
 a solution of 
\begin{equation}\label{qstrong}
	H(x,\nabla v,\nabla^2 v)-\int_{{\bf R^N}} [v(x_1+z,x_2+\gamma^{-1}z)-v(x_1,x_2)
\end{equation}
$$
	-\I \left\langle (z,\gamma^{-1}z), \nabla_{(x_1,x_2)} v(x_1,x_2) \right\rangle ] dq(z)=0\qquad x=(x_1,x_2)\in {\bf T^{2N}}, 
$$
where  $\gamma\in {\bf R}\backslash {\bf Q}$,  $H(x,p,R)$$\in$$C({\bf R^{2N}} \times{\bf R^{2N}} \times {\bf S^{2N}}  )$ satisfies 
\begin{equation}\label{positive2}
	H(x,0,O)\geq 0 \qquad \forall x\in {\bf R^{2N}},
\end{equation}
 and $dq(z)$ satisfies 
\begin{equation}\label{trans3}
	\int_{D} {\bf 1} dq(z)>0 \qquad \forall x\in {\bf R^{2N}}, \quad \forall D\subset {\bf R^N} \quad \hbox{open}.
\end{equation}

We claim the following. \\

{\bf Proposition 3.2.$\quad$} 
\begin{theorem} Let (\ref{Radon}), (\ref{positive2}), and (\ref{trans3}) hold.
Let $u$ be a periodic viscosity subsolution of  (\ref{qstrong}) in ${\bf R^{2N}}$. 
Assume 
that there exists a maximum point $\hx=(\hx_{1},\hx_{2})\in {\bf R^{2N}}$. Then, $u$ is constant almost everywhere in ${\bf R^{2N}}$. 
\end{theorem}
{\bf Proof.} We use the argument by contradiction.  Put 
$D_0=\{x\in {\bf R^{2N} }| u(x)< u(\hx)\}$. Assume that $D_0$ is non-empty and thus open. Since $\gamma$ is irrational, the set $\{\hx+(z,\gamma^{-1} z)|z\in {\bf R^N}\}$ is dense in ${\bf T^{2N}}$, and in particularly in $D_0/ [0,1]^N$.  Thus, 
$D_1$$=$$\{z\in {\bf R^N}|$$\quad  u(\hx_{1}+z,\hx_{2}+\gamma^{-1}z)$$ < u(\hx_{1},\hx_{2})  \}$ 
 is non-empty and open. On the other hand,  
since $\hx$ is a maximum point of $u$, $(0,O)\in J^{2,+}_{{\bf R^{2N}}}u(\hx)$, i.e. for any $\delta>0$ there exists $\nu>0$ such that 
$$
	u(\hx+w)\leq u(\hx) + \left\langle 0,w  \right\rangle + \frac{1}{2} \left\langle Ow,w \right\rangle +\delta |w|^2 \quad \forall |w|\leq \nu. 
$$
Take $\nu_0>0$ such that $|(z,\gamma^{-1} z)|\leq \nu$ for any $|z|<\nu_0$. 
From the definition of the viscosity subsolution, 
$$
	H(\hx,0,O) 
	- \int_{|z|<\nu_0}  \frac{1}{2}\left\langle (O+2\delta I)(z,\gamma^{-1} z),  (z,\gamma^{-1} z) \right\rangle dq(z) 
	\qquad\qquad\qquad\qquad
$$
$$
	\qquad\qquad - \int_{|z|\geq \nu_0} [ u(\hx+(z,\gamma^{-1} z))-u(\hx)
	-\I \left\langle 0,(z,\gamma^{-1} z)  \right\rangle ]dq(z) \leq 0. 
$$
Now, (\ref{positive2}), (\ref{trans3}), the definition of $D_1$, and the above inequality lead 
$$
	0< -\int_{D_1\cap \{|z|\geq \nu_0 \}} [u(\hx_1+z,\hx_2+\gamma^{-1}z)-u(\hx_1,\hx_2)]dq(z)
$$
$$
	\leq
	 -\int_{|z|\geq \nu_0} [u(\hx_1+z,\hx_2+\gamma^{-1}z)-u(\hx_1,\hx_2)
	-\I \left\langle (z,\gamma^{-1}z), 0 \right\rangle ]dq(z)
	\leq O(\delta), 
$$
which is a contradiction for $\delta>0$ sufficiently small. Therefore, $D_0$ must be  measure zero. \\

{\bf Remark 3.2.} Proposition 3.2 can be generalized to the operator in ${\bf R^{MN}}$. \\

\section{Ergodic problems for the quasi-periodic homogenizations.}

$\qquad$First, we study the following general ergodic problem, including the cases  I and II as special examples. \\

(P) Find a unique constant $d$ such that the following problem has at least a viscosity solution $v(y)$, bounded in  ${\bf R^N}$: 
$$
	d +\sup_{\alpha\in \cal{A}}\{\la -\beta(\a),\nabla v(y) \ra -c(\a) \} -a(y) \int_{{\bf R^N}} [v(y+z)-v(y)
$$
\begin{equation}\label{p1}
	- {\bf 1}_{|z|\leq 1} \left\langle z, \nabla v(y) \right\rangle] dq(z)
	-a(y)I-g_M(\gamma_1^{-1}y,...,\gamma_M^{-1}y)=0 \quad y\in {\bf R^N},
\end{equation}
where $\beta(\a)\in {\bf R^N}$ $(\a\in \cal{A})$, and 
\begin{equation}\label{c}
	|c(\a)|\leq C \qquad \forall \a\in \cal{A}. 
\end{equation}
In some cases, the number $d$ can only be characterized by the following (see \cite{al}):  for any $\mu>0$ there exist $\uv$ and $\ov$ such that 
$$
	d+\sup_{\alpha\in \cal{A}}\{\la -\beta(\a),\nabla \uv(y)\ra -c(\a)\}-a(y) \int_{{\bf R^N}} [\uv(y+z)-\uv(y)
$$
$$
	\qquad\qquad- {\bf 1}_{|z|\leq 1} \left\langle z, \nabla \uv (y) \right\rangle] dq(z)
	-a(y)I -g_M(\gamma_1^{-1}y,...,\gamma_M^{-1}y)\leq \mu\quad y\in {\bf R^N},  
$$
\begin{equation}\label{pw}\quad
\hbox{and} \qquad\qquad\qquad\qquad\qquad\qquad\qquad\qquad\qquad\qquad\qquad\qquad\qquad\qquad\qquad\qquad\qquad
\end{equation}
$$
	d+\sup_{\alpha\in \cal{A}}\{\la -\beta(\a),\nabla \ov(y)\ra -c(\a)\}-a(y) \int_{{\bf R^N}} [\ov(y+z)-\ov(y)
$$
$$
	\qquad\qquad- {\bf 1}_{|z|\leq 1} \left\langle z, \nabla \ov(y) \right\rangle] dq(z)
	-a(y)I -g_M(\gamma_1^{-1}y,...,\gamma_M^{-1}y)\geq -\mu\quad y\in {\bf R^N}. 
$$

We need also  the following formulation. \\

(Q) Find a unique constant $d$ such that the following problem has at least a viscosity solution $w(\oy)$ ($\oy=(y_1,...,y_M)$), periodic in ${\bf T^{MN}}$, 
$$
	d+\sup_{\alpha\in \cal{A}}\{\la -B(\a),\nabla w(\oy)\ra -c(\a) \}-a(y_1)
	\int_{{\bf R^N}} [
	w(\oy+\Gamma^{-1}z)-w(\oy)
$$
$$
	\qquad-\I \left\langle \Gamma^{-1}z, \nabla w(\oy) \right\rangle ]dq(z) 
	-a(y)I -g_M(y_1,...,y_M)=0\qquad \oy\in {\bf T^{MN}}, 
$$
where $B(\a)=(\gamma_{1}^{-1}\beta(\a),...,	\gamma_{M}^{-1}\beta(\a))$, $\Gamma^{-1}=(\gamma_{1}^{-1},...,	\gamma_{M}^{-1})$. \\

 In fact, (P) and (Q) are related by $v(y)$$=w(\gamma_{1}^{-1} y,...,\gamma_M^{-1}y)$. We abbreviate the weaker version of (Q). \\
	
	As  in \cite{al}, we approximate  (P) by 
$$
	\l v_{\l} +\sup_{\alpha\in \cal{A}}\{\la -\beta(\a),\nabla v_{\l}(y) \ra -c(\a) \} -a(y) \int_{{\bf R^N}} [v_{\l}(y+z)-v_{\l}(y)\qquad\qquad\qquad\qquad
$$  
\begin{equation}\label{vl}
	- {\bf 1}_{|z|\leq 1} \left\langle z, \nabla v_{\l}(y) 	 		       
    \right\rangle] dq(z)-a(y)I-
	g_M(\gamma_1^{-1}y,...,\gamma_M^{-1}y)=0 \quad y\in {\bf R^N}, 
\end{equation}
and (Q) by 
$$
	\l w_{\l}+\sup_{\alpha\in \cal{A}}\{\la -B(\a),\nabla w_{\l}(\oy)\ra -c(\a) \}-a(y_1)
	\int_{{\bf R^N}} [
	w_{\l}(\oy+\Gamma^{-1}z)-w_{\l}(\oy)
$$
\begin{equation}\label{wl}
	\qquad-\I \left\langle \Gamma^{-1}z, \nabla w_{\l}(\oy) \right\rangle ]dq(z) 
	-a(y)I -g_M(\oy) =0\qquad \oy\in {\bf T^{MN}},
\end{equation}
 for $\l\in (0,1)$. \\

{\bf Remark 4.1.} If $dq(z)=\frac{1}{|z|^{N+1}}dz$ ($\a=1$ in (\ref{vl})), then the L{\'e}vy operator is close to the first-order partial differential opertor. 
Certainly, a condition is necessary between $(\beta,c)$ and $a(\cdot)$ to determine which term: $\sup_{\alpha\in \cal{A}}\{\la -B(\a),\nabla w_{\l}(\oy)\ra -c(\a) \}$, and 
$-a(y_1)
	\int_{{\bf R^N}} [
	w_{\l}(\oy+\Gamma^{-1}z)-w_{\l}(\oy) -\I \left\langle \Gamma^{-1}z, \nabla w_{\l}(\oy) \right\rangle ]dq(z)$, in major, serves for the ergodicity. This is not a trivial question, 
and to avoid the complexity we assume that $a(\cdot)\equiv a$, if $\a=1$. \\

Our claim is the following. \\

{\bf Theorem 4.1.$\quad$} 
\begin{theorem} Assume that (\ref{a}), (\ref{nonreso}), (\ref{holder}), (\ref{gholder}), (\ref{Radon}), (\ref{trans}), and (\ref{c}) hold.  
 Assume also that either $a(\cdot)\equiv a$ ($a>0$ is a constant), or $\beta(\a) \equiv 0$ ($\forall \a \in {\mathcal A}$). 
 Let $v_{\l}$ be the solution of (\ref{vl}). 
	Then, for any $\theta\in (0,\theta_0]$, there exists a constant $C_{\theta}>0$ independent on $\l>0$, such that
\begin{equation}\label{hcont}
	|v_{\l}(y)-v_{\l}(y')|\leq \frac{C_{\theta}}{\l} |y-y'|^{\theta} \qquad \forall y,y'\in {\bf R^N}. 
\end{equation}
\end{theorem}

We prepare some lemmas. \\

{\bf Lemma 4.2.$\quad$}  
\begin{theorem} Consider (\ref{vl}) (resp. (\ref{wl})). The following hold. \\
(i) Let $\uv_{\l}$ (resp. $\uw_{\l}$) be a bounded USC subsolution of (\ref{vl}) (resp. (\ref{wl})). 
Let $\ov_{\l}$ (resp. $\ow_{\l}$) be a bounded LSC supersolution of (\ref{vl}) (resp. (\ref{wl})). 
Then, $\uv_{\l}\leq \ov_{\l}$ (resp. $\uw_{\l}\leq \ow_{\l}$) holds in ${\bf R^N}$ (resp. ${\bf T^{MN}}$). \\
(ii) There exists a unique bounded viscosity solution $v_{\l}$ (resp. $w_{\l}$) of  (\ref{vl}) (resp. (\ref{wl})). \\
\end{theorem}

{\bf Proof.} (i) The proof of the comparison principle can be done by a standard way. We refer the readres to \cite{ar3}, \cite{ar5}, \cite{bbp}, 
 \cite{bi}. \\
(ii) The existence of the solutions can be shown by the Perron's method (see Crandall, Ishii and Lions \cite{users}).  We refer the readres to \cite{ar3}, \cite{ar5}, \cite{bbp}, \cite{bi} for details.  \\

We multiply (\ref{vl}) by $\l>0$, put $m_{\l}=\l v_{\l}$, $f(y)=g_M(\gamma_1^{-1}y,...,\gamma_M^{-1}y)$$+a(y)I$ to have 
\begin{equation}\label{ml}
	\l m_{\l} (y)+\sup_{\alpha\in \cal{A}}\{\la -\beta(\a),\nabla m_{\l}(y) \ra -\l c(\a)\}
	 -a(y)\int_{{\bf R^N}} [m_{\l}(y+z)-m_{\l}(y)
\end{equation}
$$
	\qquad\qquad\qquad\qquad\qquad
	-\I \left\langle z, \nabla m_{\l}(y) \right\rangle ]dq(z)= \l f(y) \qquad y\in {\bf R^N}. 
$$

{\bf Lemma 4.3.$\quad$} 
\begin{theorem} There exists a constant $M>0$ such that 
\begin{equation}\label{vM}
	|m_{\l}|\leq M\qquad \forall \lambda\in (0,1), 
\end{equation}
and  for any $\theta\in (0,\theta_0]$ there exists a constant $C_{\theta}>0$ independent on $\l>0$ such that
\begin{equation}\label{hcont}
	|m_{\l}(y)-m_{\l}(y')|\leq C_{\theta} |y-y'|^{\theta} \qquad \forall y,y'\in {\bf R^N}. 
\end{equation}
\end{theorem}

{\bf Proof.} Since $m_{\l}=\l v_{\l}$, from the comparison principle for (\ref{vl}) (Lemma 4.2 (i)),  (\ref{vM}) is clear. 
Fix $\theta\in (0,\theta_0]$, and let $r_0>0$ be a constant to be determined later. Put 
\begin{equation}\label{Ctheta}
	C_{\theta}=\frac{2M}{r_0^{\theta}}. 
\end{equation}
We use the argument by contradiction to prove the claim. So, assume that there exist $\tilde{y}$, $\ty'$ such that 
\begin{equation}\label{contra}
	m_{\l}(\ty)-m_{\l}(\ty')> C_{\theta} |\ty-\ty'|^{\theta}. 
\end{equation}
From (\ref{vM}), (\ref{Ctheta}), 
 $|\ty-\ty'|<r_0$ must hold. We regularize $m_{\l}$ by the sup-convolution and the  inf-convolution: for $\l>0$
$$
	m^r(y)=\sup_{y'\in {\bf R^N}} \{ m_{\l}(y')-\frac{r}{2}|y-y'|^2\},\quad
	m_r(y)=\inf_{y'\in {\bf R^N}} \{ m_{\l}(y')+\frac{r}{2}|y-y'|^2\}. 
$$
Remark that (see \cite{ar5}) for any $\nu>0$ we can take $r>0$ small enough so that 
$$
	\l m^{r}+\sup_{\alpha\in \cal{A}}\{\la -\beta(\a),\nabla m^{r}(y) \ra -\l c(\a) \}
	-a(y)\int_{{\bf R^N}} [m^{r}(y+z)-m^{r}(y)
$$
$$
	\qquad\qquad\qquad\qquad\qquad\qquad\qquad\qquad\qquad   
	- \I \left\langle z, \nabla m^{r}(y) \right\rangle ]dq(z)\leq \l  f(y) + \nu,
$$
$$
	\l m_{r}+\sup_{\alpha\in \cal{A}}\{\la -\beta(\a),\nabla m_{r}(y) \ra -\l c(\a)\}
	-a(y)\int_{{\bf R^N}} [m_{r}(y+z)-m_{r}(y)
$$
$$
	\qquad\qquad\qquad\qquad\qquad\qquad\qquad\qquad\qquad\qquad
	- \I \left\langle z, \nabla m_{r}(y) \right\rangle ]dq(z)\geq \l  f(y) - \nu, 
$$
in the sense of the viscosity solution. Since $m_{r}\leq m_{\l} \leq m^{r}$,  from (\ref{contra})
\begin{equation}\label{mlcont}
	m^{r}(\ty)- m_{r}(\ty') > C_{\theta}|\ty-\ty'|^{\theta}. 
\end{equation}
Define 
$$
	\Phi(y,y')=m^{r}(y)- m_{r}(y') - C_{\theta}|y-y'|^{\theta} 
	\quad \forall (y,y') \in {\bf R^{2N}}. 
$$
  Since $m^r$, $m_r$ are the sup and the inf convolutions of $\l v_{\l}= \l w_{\l}(\gamma_1^{-1}y,...,\gamma_M^{-1}y)$, and since $w_{\l}$ is periodic,  the maximum point of $\Phi$ exists. Let $(\hy,\hy')$ be the maximum point of $\Phi$.

  Put  $p=\nabla_y \phi(\hy,\hy')$, $Q=\nabla^2_y \phi(\hy,\hy')$. 
 In particular, we may assume that $(\hy,\hy')$ is a global strict maximum point of $\Phi$.   We can take an open precompact subset $\mathcal{O}\subset {\bf R^{2N}}$ such that $(\hy,\hy')\in \mathcal{O}$, $\sup_{\mathcal{O}}\Phi(y,y')-\sup_{\p \mathcal{O}} \Phi(y,y')>0$. 
Then, from the Alexandrov's maximum principle and the Jensen's lemma (see Fleming and Soner \cite{fs}),  the following holds (\cite{ar5}). \\

{\bf Lemma A.$\quad$(\cite{ar5}$\quad$ Lemma 1.3.)} 
\begin{theorem} 
(i) There exists a sequence $(y_j,y'_j)$ in $\mathcal{O}$ which converges  to 
 $(\hy,\hy')$ as $j\to \infty$, and $(p_j,Y_j)\in J^{2+}_{\Omega} m^r(y_j)$,  $(p'_j,Y'_j)\in J^{2-}_{\Omega} m_r(y'_j)$  such that 
$$
	\lim_{j\to \infty} p_j=\lim_{j \to \infty} p_j'=p,\quad Y_{j}\leq Y'_{j} \quad \forall j\in {\bf Z}. 
$$
(ii) For $P_j=(p_j-p,-(p'_j-p))$, $\Phi_j(y,y')=\Phi(y,y')-\la P_j,(y,y') \ra$ takes a maximum at $(y_j,y'_j)$ in $\mathcal{O}$. \\
(iii) For any $z\in {\bf R^N}$ such that $(y_j+z,y'_j+z)\in \mathcal{O}$
$$
	 m^r(y_j+z)- m^r(y_j)-\la p_j,z \ra \leq  m_r(y'_j+z)- m_r(y'_j)-\la p'_j,z \ra. 
$$
\end{theorem}

Take a pair of positive numbers $(\nu_j,\delta_j)$ such that 
$$
	m^{r}(y_j+z)\leq  m^{r}(y_j) + \left\langle z, p_j \right\rangle + \frac{1}{2}\left\langle Y_j z,z \right\rangle  + \delta_j |z|^2
	\qquad \forall |z|\leq \nu_j, 
$$
$$
	m_{r}(y'_j+z)\geq  m_{r}(y'_j) + \left\langle z, p'_j \right\rangle + \frac{1}{2}\left\langle Y'_j z,z \right\rangle  - \delta_j |z|^2
	\qquad \forall |z|\leq \nu_j. 
$$
From the definition of the viscosity solution,  by remarking that $a(y_j)$, $a(y'_j)$$\geq a_0>0$, 
$$
	\frac{\l m^{r}(y_j)}{a(y_j)}+\sup_{\alpha\in \cal{A}}\{\la -\frac{\beta(\a)}{a(y_j)}, p_j \ra -\frac{\l c(\a)}{a(y_j)} \}  
	 -\int_{|z|<\nu_j} \frac{1}{2}\left\langle (Y_j+2\delta_j I)z,z \right\rangle dq(z)
$$
$$
	-\int_{|z|>\nu_j} [m^{r}(y_j+z)-m^{r}(y_j)-\I \left\langle z, p_j \right\rangle ] dq(z) \leq  \frac{\l f(y_j) + \nu}{a(y_j)},
$$
$$
	\frac{\l m_{r}(y'_j)}{a(y_j')}+\sup_{\alpha\in \cal{A}}\{\la -\frac{\beta(\a)}{a(y_j')}, p'_j \ra -\frac{\l c(\a)}{a(y_j')} \}   -\int_{|z|<\nu_j} \frac{1}{2}\left\langle (Y'_j-2\delta_j I)z,z \right\rangle dq(z)
$$
$$
	-\int_{|z|>\nu_j} [m_{r}(y'_j+z)-m_{r}(y'_j)-\I \left\langle z, p'_j \right\rangle ] dq(z) \geq  \frac{\l f(y'_j) - \nu}{a(y_j')}. 
$$
We take the difference of  two inequalities. 
Put 
$$
	\mathcal{O}_j^c =  \{|z|>v_j\} \cap \{ z|\quad (y_j+z,y_j'+z)\in \mathcal{O}^c\}. 
$$
By remarking $Y_j\leq Y_j'$, Lemma A (iii), by choosing $\a\in \mathcal A$ appropriately 
$$
	\frac{\l a(y_j')m^{r}(y_j) - \l a(y_j) m_{r}(y'_j)}{a(y_j)a(y_j')}
	-\la \frac{\beta(\a)}{a(y_j)}, p_j \ra + \la \frac{\beta(\a)}{a(y_j')}, p_j' \ra
$$
$$
	\leq 
	- \int_{ z\in \mathcal{O}_j^c }
	 [m^{r}(y_j+z)-m^{r}(y_j) - m_{r}(y'_j+z)+m_{r}(y'_j)] dq(z)
$$
$$
	+ \frac{\l c(\a) |a(y_j)-a(y_j')| }{a(y_j)a(y_j')}+ \frac{\l(a(y_j')f(y_j)-a(y_j)f(y_j')) + \nu(a(y_j')-a(y_j)) }{a(y_j)a(y_j')}. 
$$
We may assume $\nu_j\to 0$ as $j\to \infty$. Then, $\mathcal{O}_j^c\to \mathcal{O}^c$.  
We let $j\to \infty$ in the above inequality, by remarking 
$$
	m^{r}(\hy)-m_{r}(\hy') -C_{\theta}|\hy-\hy'|^{\theta}\geq  m^{r}(\hy)-m_{r}(\hy')-C_{\theta}|\hy-\hy'|^{\theta}, 
$$
and by multiplying  by $a(\hy)a(\hy')$ 
$$
	\l (a(\hy')m^{r}(\hy)-a(\hy)m_{r}(\hy')) + (a(\hy)-a(\hy')) \left\langle \beta(\a),p \right\rangle 
$$
$$
	\leq \l c(\a) |a(\hy)-a(\hy')|  
	+ \l(a(\hy')f(\hy)-a(\hy)f(\hy')) + \nu(a(\hy')-a(\hy)). 
$$
Since either $a(\cdot)\equiv a$,  or $\beta(\a) \equiv 0$ ($\forall \a \in {\mathcal A})$, and since  
$\nu>0$ is arbitrary, deviding the both hands side by $\l>0$, 
$$
	a_0(m^{r}(\hy)-m_{r}(\hy'))\leq C ( |a(\hy')-a(\hy)| +  |f(\hy')-f(\hy)| ). 
$$
From the H{\"{o}}lder continuity of $a$, $f$ and (\ref{mlcont}), we get 
$$
	C_{\theta}|\hy-\hy'|^{\theta} \leq C|\hy-\hy'|^{\theta_0}. 
$$
From (\ref{Ctheta}),
$\frac{2M}{r_0^{\theta_0}}\leq C |\hx-\hy|^{\theta_0-\theta}$, that is 
$$
	2M\leq C|\hx-\hy|^{\theta_0-\theta} r_0^{\theta_0}\leq Cr_0^{\theta_0}.
$$ 
Therefore, for $r_0>0$ small enough such that $r_0^{\theta_0}<\frac{M}{C}$, we get a contradiction. For such $r_0>0$, by defining 
 $C_{\theta}$ as in (\ref{Ctheta}),  we proved our claim. \\

{\bf Proof of Theorem 4.1.} Since $m_{\l}=\l v_{\l}$ in Lemma 4.2, the claim of Theorem 4.1 is clear. \\

{\bf Theorem 4.4.$\quad$} 
\begin{theorem} Assume that (\ref{a}), (\ref{nonreso}), (\ref{holder}),  (\ref{gholder}), (\ref{Radon}), (\ref{trans}), and (\ref{c}) hold.  
 Assume also that either $a(\cdot)\equiv a$ ($a>0$ is a constant), or $\beta(\a) \equiv 0$ ($\forall \a \in {\mathcal A}$).  
Let $v_{\l}$ be the solution of (\ref{vl}).  
Then, there exists a unique number $d$ such that 
$$
	\lim_{\l \to 0} \l v_{\l} (y)=d \quad \hbox{uniformly in} \quad {\bf R^N}, 
$$
which  is characterized by (\ref{pw}). 
\end{theorem}

{\bf Proof.} From Theorem 4.1, by the Ascoli-Arzela theorem, there exists a subsequence  $\l'\to 0$ such that
$$
	\lim_{\l' \to 0} \l' v_{\l'} (y)=d(y) \qquad y\in {\bf R^N}. 
$$
We still use $\l$ instead of $\l'$ to simplify the notation. Remark that $d(y)$ is  H{\"{o}}lder continuous. 
Since $v_{\l}(y)$$=w_{\l}(\gamma_1^{-1}y,...,\gamma_M^{-1}y)$ and $w_{\l}$ is periodic, 
the above convergence  is uniform in ${\bf R^N}$.   
Multiplying (\ref{vl}) by $\l>0$, passing $\l\to 0$, we get 
$$
	\sup_{\alpha\in \cal{A}}\{\la -\beta(y,\a), \nabla d(y) \ra \}   - a(y) \int_{{\bf R^N}} [d(y+z)-d(y)-\I \left\langle z,\nabla d(y) \right\rangle ]dq(z) =0. 
$$
Since $d(y)$ is a uniform  limit of a sequence of quasi-periodic functions $\l v_{\l}(y)$, it takes a maximum at some point $\hy \in {\bf R^N}$.  From the strong maximum principle in 
Theorem 3.1,  $d(y)\equiv d$. 
The uniqueness of $d$ can be proved by the standard argument (see \cite{al} for example). Let $\mu>0$ be arbitrary. 
From the uniform convergence of $\l v_{\l}$ as $\l$ goes to $0$,  for $\l>0$ small enough if we put $\uv=v_{\l}$ and $\ov=v_{\l}$ they satisfy 
 (\ref{pw}). The claims in Theorem 4.4 are thus proved. \\

Next, we study the ergodic problem of the first-order PDE, which includes the case I as a special example.  
Let us consider the following deterministic system 
\begin{equation}\label{dynamicsb}
	\frac{d y_{\a}}{dt}=\beta(\a(t)) \quad t>0, \quad y_{\a}(0)=y\in {\bf T^N}, 
\end{equation}
where $\a(\cdot)$ is a measurable function from $[0,\infty)$ to $\mathcal A$, which we call a control. We assume the following. \\

(A) A controlled dynamical system (\ref{dynamicsb}) is approximately controllable if for any $y$, $y'\in {\bf T^N}$, and for any $\delta>0$, there exists a control 
$\a(\cdot)$ and $T_{\delta}>0$ such that the solution $y_{\a}(t)$ of (\ref{dynamicsb}) satisfies $|y'-y(T_{\delta})|<\delta$. \\

Under the above condition, we intend to solve the following. \\

(R) Find a unique constant $d$ such that the following problem has at least a viscosity solution $v(y)$, bounded in  ${\bf R^N}$: 
$$
	d +\sup_{\alpha\in \cal{A}}\{\la -\beta(\a),\nabla v(y) \ra -c(\a)\} -a(y_1)I -g_M(\gamma_1^{-1}y,...,\gamma_M^{-1}y)=0 \quad y\in {\bf R^N}. 
$$

(S) Find a unique constant $d$ such that the following problem has at least a viscosity solution $w(\oy)$ ($\oy=(y_1,...,y_M)$), periodic in ${\bf T^{MN}}$: 
$$
	d+\sup_{\alpha\in \cal{A}}\{\la -B(\oy,\a),\nabla w(\oy)\ra -c(\a) \}-a(y_1)I 
	-g_M(\oy) =0\qquad \oy\in {\bf T^{MN}}, 
$$
where $B(\a)=(\gamma_{1}^{-1}\beta(\a),...,	\gamma_{M}^{-1}\beta(\a))$, $\Gamma^{-1}=(\gamma_{1}^{-1},...,	\gamma_{M}^{-1})$. \\

 The problems (R) and (S) are related by $v(y)$$=w(\gamma_{1}^{-1} y,...,\gamma_M^{-1}y)$. We abbreviate the weaker versions of (R) and (S).  
As before, we approximate the problems, for $\l\in (0,1)$ \\
$$
	 \l v_{\l}+\sup_{\alpha\in \cal{A}}\{\la -\beta(\a),\nabla v_{\l}(y) \ra-c(\a) \} -a(y_1)I -g_M(\gamma_1^{-1}y,...,\gamma_M^{-1}y)=0 \qquad
$$
\begin{equation}\label{approxR}
\qquad\qquad\qquad\qquad\qquad\qquad\qquad\qquad\qquad\qquad\qquad\qquad
	\quad y\in {\bf R^N}. 
\end{equation}
$$
	\l w_{\l} +\sup_{\alpha\in \cal{A}}\{\la -B(\oy,\a),\nabla w_{\l}(\oy)\ra -c(\a) \}-a(y_1)I 
	-g_M(\oy)=0,\qquad \qquad\quad
$$
\begin{equation}\label{approxS}\qquad\qquad\qquad\qquad\qquad\qquad\qquad\qquad\qquad\qquad\qquad\qquad
	\quad y\in {\bf T^{MN}}. 
\end{equation}
{\bf Theorem 4.5.$\quad$} 
\begin{theorem} Assume that (\ref{nonreso}), (\ref{holder}), (\ref{gholder}), and (\ref{c}) hold. Assume also that (\ref{dynamicsb}) is approximately controllable. 
Let $v_{\l}$ be the solution of (\ref{approxR}).  
Then, there exists a unique number $d$ such that 
$$
	\lim_{\l \to 0} \l v_{\l} (y)=d \quad \hbox{uniformly in} \quad {\bf R^N}. 
$$
Moreover, the number $d$ is characterized by  (\ref{pw}). 
\end{theorem}

{\bf Proof.}  Let $w_{\l}$ be the periodic solution of (\ref{approxS}). \\
 (Step 1)  We first show that there exists a constant $C>0$ such that
\begin{equation}\label{wcont}
	|\l w_{\l}(\oy)-\l w_{\l}(\oy')|\leq C |\oy-\oy'|^{\theta_0}\quad \forall \oy,\oy'\in {\bf T^{MN}}, \quad \forall \l\in (0,1). 
\end{equation}
Let $\a(t)$ be an arbitrary measurable function from $(0,\infty)$ to ${\mathcal A}$. Let $Y_{\a}(t)$,  $Y'_{\a}(t)$ be respectively the solution of 
$$
	\frac{d Y_{\a}}{dt}=B(\a(t)) \quad t>0, \quad Y_{\a}(0)=\oy; \quad
	\frac{d Y'_{\a}}{dt}=B(\a(t)) \quad t>0, \quad Y_{\a}(0)=\oy'. 
$$
Remark that there exists a constant $L>0$ such that 
\begin{equation}\label{L}
	|Y_{\a}(t)-Y'_{\a}(t)|\leq L|\oy-\oy'|\quad \forall t\geq 0. 
\end{equation} 
Put $f(\oy)=a(y_1)I $$+g_M(\oy)$. Since (see for example \cite{users})
$$
	 w_{\l}(\oy)=\inf_{\a(\cdot)} \{
	\int_0^{\infty} e^{-\l t} (f(Y_{\a}(t))+c(\a(t)) )dt
	\},
$$
$$
	 w_{\l}(\oy')=\inf_{\a(\cdot)} \{
	\int_0^{\infty} e^{-\l t} (f(Y'_{\a}(t))+c(\a(t))) dt
	\},
$$
for any $\nu>0$, we can take a control $\a(\cdot)$ such that 
$$
	 w_{\l}(\oy)-w_{\l}(\oy') \leq 
	\int_0^{\infty} e^{-\l t} |f(Y_{\a}(t)) - f(Y'_{\a}(t))| dt + \nu
	\leq \frac{L}{\l} |\oy-\oy'|^{\theta_0} +\nu,
$$
where we used (\ref{holder}), (\ref{gholder}), (\ref{L}) to derive the last inequality. Since $\nu>0$ is arbitrary, (\ref{wcont}) is shown. \\
(Step 2) From (\ref{wcont}), we can extract a subsequence $\l'\to 0$ such that 
$$
	\lim_{\l'\to 0} \l'  w_{\l '}(\oy) =d(\oy) \quad \hbox{uniformly in}\quad  \oy\in {\bf T^{MN}},
$$
where $d(\oy)$ is H\"{o}lder continuous and periodic. Multiplying (\ref{approxS}) by $\l'>0$, and tending $\l'$ to zero, we deduce that $d(\oy)$ satisfies 
\begin{equation}\label{limit}
	\sup_{\a\in \mathcal A} \{
	\left\langle -B(\a), \nabla d(\oy)\right\rangle 
	\} =0 \quad \oy\in {\bf T^{MN}}. 
\end{equation}
Now, 
$$
	\left\langle -B(\a), \nabla d(\oy)\right\rangle  \leq 0 \quad \forall \a\in {\mathcal A}. 
$$ 
Since $d(\cdot)$ is Lipshitz,  the above holds almost everywhere in ${\bf T^{MN}}$.  Then, since  $d(\cdot)$ is periodic, 
$$
	\left\langle -B(\a), \nabla d(\oy)\right\rangle  = 0 \quad  \forall \a\in   {\mathcal A}. 
$$
Thus, 
$$
	d(\oy+\int_0^t B(\a(s)) ds)=d(\oy) \qquad \forall t\geq 0, \quad \forall  \a(\cdot):[0,\infty)\to {\mathcal A}. 
$$
Remark that $B(\a)=(\gamma_{1}^{-1}\beta(\a),...,	\gamma_{M}^{-1}\beta(\a))$, and the set $\{\gamma_i\}$ ($1\leq i\leq M$) satisfies the 
non-resonance condition.  From (A), the set 
$$
	\bigcup_{\a(\cdot),t\geq 0} \{y_1+\int_0^t \beta (\a(s)) ds \} \quad \hbox{is dense in}\quad \bf T^N. 
$$ 
Thus, the set
$$
	\bigcup_{\a(\cdot),t\geq 0}\{\oy+  \int_0^t B(\a(s)) ds \} \quad \hbox{is dense in}\quad \bf T^{MN}. 
$$ 
Therefore, $d(\oy)\equiv d$ in $\oy\in {\bf T^{MN}}$. \\
(Step 3) The number $d$ is unique, and it does not depend on the choice of the subsequence $\l'\to 0$.  
The proof is standard, and we refer the readres to \cite{al}. By defining $v_{\l}(y)$$=w_{\l}(y,\gamma_1^{-1}y,...,\gamma_M^{-1}y)$, we have shown the claim. \\

{\bf Proposition 4.6.$\quad$} 
\begin{theorem} Assume that (\ref{a}), (\ref{nonreso}), (\ref{holder}), and (\ref{gholder}) hold.  
Let $(x,p,I)$$\in {\bf \Omega \times R^N \times R}$ be arbitrarily fixed. 
Then, the following hold. \\

(i) Let $\a\in (0,1)$ and $b(x,\a)\not \equiv 0$. Let (\ref{dynamicsb}) with $\beta(\a)$$=b(x,\a)$  satisfy (A). There exists a unique number $d_{x,p,I}$ which satisfies  (\ref{cell1}) in the  sense of (\ref{pw}). \\
(ii) Let $\a=1$, $b(x,\a)\not \equiv 0$, and $a(y)\equiv a$ ($y\in {\bf R^N}$). There exists a unique number $d_{x,p,I}$ which satisfies (\ref{cell2}) in the sense of (\ref{pw}). \\
(iii) Let  $\a\in (1,2)$, or $\a\in (0,1]$ and $b(x,\a)\equiv 0$.   There exists a unique number $d_{x,p,I}$ which satisfies   (\ref{cell3}) in the sense of (\ref{pw}). \\
\end{theorem}

{\bf Proof.} (i) Put $\beta(\a)$$=b(x,\a)$ in (\ref{approxR}), $c(\a)=\left\langle b(x,\a),p \right\rangle$. Then, from Theorem 4.5, the statement follows. \\
(ii) In (\ref{vl}), put $\beta(\a)=b(x,\a)$, $c(\a)=\left\langle b(x,\a),p \right\rangle $. The claim follows from Theorem 4.4. \\
(iii) In (\ref{vl}), put $\beta(\a)=0$, $c(\a)=0$. From Theorem 4.4, the claim follows. \\

\section{Ergodic problems for the almost periodic homogenizations. }

$\quad$Next, we solve the ergodic cell problem for  the almost periodic homogenizations. Similar to  the case of the quasi-periodic homogenizations, 
 the situation is devided into the following. 
\begin{itemize}
\item I'. $\a\in (0,1)$ and $b(x,\a)\not \equiv 0$.
\item II'. $\a=1$ and $b(x,\a)\not \equiv 0$.
\item III'. $\a\in (1,2)$, or $\a\in (0,1]$ and $b(x,\a)\equiv 0$.
\end{itemize}
The formal asymptotic expansion (see \S 2) leads, for the case I' 
\begin{equation}\label{acell1}
	d_{x,p,I} +\sup_{\alpha\in \cal{A}}\{\la -b(x,\a),p+\nabla_y v(y) \ra \}- a(y) I -g(y)=0 \quad y\in {\bf R^N}. 
\end{equation}
For the case II', 
\begin{equation}\label{acell2}
	d_{x,p,I} 
	+\sup_{\alpha\in \cal{A}}\{\la -b(x,\a),p+\nabla_y v(y)\ra \}-a(y) \int_{{\bf R^N}} [v(y+z)-v(y)
\end{equation}
$$
	- {\bf 1}_{|z|\leq 1} \left\langle z, \nabla_y v(y) \right\rangle] \frac{1}{|z|^{N+\a}}dz -a(y)I-g(y)=0 \quad y\in {\bf R^N}. 
$$
And for the case III', 
\begin{equation}\label{acell3}
	d_{x,p,I}- a(y)\int_{{\bf R^N}} [v(y+z)-v(y)
	- \I \left\langle z, \nabla_y v(y) \right\rangle] \frac{1}{|z|^{N+\a}}dz- a(y)I 
\end{equation}
$$
	\qquad\qquad\qquad\qquad\qquad\qquad\qquad\qquad \qquad - g(y)=0 \qquad y\in {\bf R^N}. 
$$

Our claim is the following. \\

{\bf Proposition 5.1.$\quad$} 
\begin{theorem} Assume that (\ref{a}), (\ref{holder}), and (\ref{gholder}) hold, and that $g$ in (\ref{almostex}) is uniformly almost periodic in the sense of Bohr 
 in ${\bf R^N}$. 
Let $(x,p,I)$$\in \Omega \times {\bf R^N} \times {\bf R}$ be arbitrarily fixed. 
Then, the following hold. \\

(i) Let $\a\in (0,1)$ and $b(x,\a)\not \equiv 0$. Let (\ref{dynamicsb}) with $\beta(\a)$$=b(x,\a)$  satisfy (A). There exists a unique number $d_{x,p,I}$ which satisfies  (\ref{acell1}) in the  sense of (\ref{pw}). \\
(ii) Let $\a=1$, $b(x,\a)\not \equiv 0$, and $a(y)\equiv a$ ($y\in {\bf R^N}$). There exists a unique number $d_{x,p,I}$ which satisfies (\ref{acell2}) in the sense of (\ref{pw}). \\
(iii) Let $\a\in (1,2)$, or $\a\in (0,1]$ and $b(x,\a)\equiv 0$. There exists a unique number $d_{x,p,I}$ which satisfies   (\ref{acell3}) in the sense of (\ref{pw}). \\
\end{theorem}

 We consider the following general ergodic problem, which includes (\ref{acell2}) and (\ref{acell3}) as special cases. Find a unique constant 
 $d>0$ such that there exists a bounded  viscosity solution $v$ of 
\begin{equation}\label{ae1}
	d+\sup_{\a\in \cal{A}} \{ \la \beta(\a),  \nabla_y v (y) \ra -c(\a) \}
	- a(y) \int_{{\bf R^N}} [v(y+z)-v(y)\qquad\qquad\qquad
\end{equation}
$$
	\qquad\quad\qquad
	-\I \left\langle z,\nabla_y v(y) \right\rangle ] dq(z)-a(y)I- g(y) =0 \qquad y\in {\bf R^N}. 
$$

For (\ref{acell1}), we are interested in finding a unique constant $d>0$ such that there exists a bounded  viscosity solution $v$ of 
\begin{equation}\label{ae2}
	d+\sup_{\a\in \cal{A}} \{ \la \beta(\a),  \nabla_y v (y) \ra -c(\a) \}-a(y)I- g(y) =0 \qquad y\in {\bf R^N}. 
\end{equation}
 
In some cases,  the number $d$ satisfies (\ref{ae1}) (resp.(\ref{ae2})) in the sense of (\ref{pw}).  We use a useful characterization of the uniformly almost periodic function in  Braides \cite{braides}, which was first shown by Bohr (see \cite{besicovitch}) for the one dimensional case. 

\begin{theorem} {\bf Lemma B. (\cite{braides} Definition A.1, Theorem A.6.)}  If a continuous function $f(x)$ defined in ${\bf R^N}$ is uniformly almost periodic in the sense of Bohr, then $f$ 
is the uniform limit of a sequence of trigonometric polynomials. The converse is also true. 
\end{theorem}

	We refer the readers to \cite{besicovitch} and \cite{braides} for details. \\

{\bf Lemma 5.2.$\quad$} 
\begin{theorem} 
(i) Assume that (\ref{a}), (\ref{holder}),  (\ref{gholder}), (\ref{Radon}), (\ref{trans}), and (\ref{c}) hold.   
Assume also that  $g$ is uniformly almost periodic in ${\bf R^N}$. 
There exists a unique constant $d$ which satisfies (\ref{ae1}) in the sense of (\ref{pw}). \\
(ii) Assume that  (\ref{holder}), (\ref{gholder}), and (\ref{c}) hold, that (\ref{dynamicsb}) is approximately controllable.  
Assume also that  $g$ is uniformly almost periodic in ${\bf R^N}$. 
There exists a unique constant $d$ which satisfies (\ref{ae2}) in the sense of (\ref{pw}).
\end{theorem}

{\bf Proof.} (i) From Lemma B,   there exist a sequence of periodic functions $g_M$ ($M=1,2,...$)  
 defined in $(y_1,...,y_M)$$\in {\bf T^{MN}}$, and a sequence of numbers $\gamma_i$ $(i\in {\bf N})$ satisfying the non-resonance condition 
(\ref{nonreso}), such that 
\begin{equation}\label{uniform}
	g(y)=\lim_{M\to \infty} g_M(\gamma_1^{-1}y,...,\gamma_M^{-1}y) \quad \hbox{uniformly in}\quad {\bf R^N}. 
\end{equation}
From Theorem 4.4, for each $M$,  there exist a constant $d_M$,  $\uv_M \in USC(\bf R^N)$,  and $\ov_M \in LSC(\bf R^N)$ which satisfy : 
$$
	d_M+\sup_{\a\in \cal{A}} \{ \la \beta(\a),  \nabla_y \uv_M (y) \ra -c(\a)\}
	- a(y) \int_{{\bf R^N}} [\uv_M(y+z)-\uv_M (y)\qquad\qquad\qquad
$$
$$
	\qquad\quad\qquad
	-\I \left\langle z,\nabla_y \uv_M (y) \right\rangle ] dq(z)-g_M(y)\leq \mu \qquad y\in {\bf R^N}, 
$$
$$
	d_M+\sup_{\a\in \cal{A}} \{ \la \beta(\a),  \nabla_y \ov_M (y) \ra -c(\a) \}
	- a(y) \int_{{\bf R^N}} [\ov_M (y+z)-\ov_M (y)\qquad\qquad\qquad
$$
$$
	\qquad\quad\qquad
	-\I \left\langle z,\nabla_y \ov_M (y) \right\rangle ] dq(z)-g_M(y)\geq -\mu \qquad y\in {\bf R^N}. 
$$

Remarking that there exists a constant $C>0$ such that 
$$
	|d_M|\leq C\qquad \forall M>0,
$$
for $g(y)$, $g_M(\Gamma^{-1} y)$ ($M\in {\bf N}$) are uniformly bounded in ${\bf R^N}$. Thus, 
we can extract a sequence $M'\to \infty$ such that $\lim_{M'\to \infty} d_{M'}=d$. 
 Define $v^{\ast}(y)=\overline{\lim}_{M'\to \infty} \uv_{M'} (y)$, $v_{\ast}(y)=\underline{\lim}_{M'\to \infty} \ov_{M'} (y)$. 
From Barles and Perthame \cite{bp1},  by passing $M'\to \infty$ in the above inequalities, we find that $v^{\ast}$ and $v_{\ast}$ respectively satisfy 
$$
	d+\sup_{\a\in \cal{A}} \{ \la \beta(\a),  \nabla_y v^{\ast} (y) \ra -c(\a)\}
	- a(y) \int_{{\bf R^N}} [v^{\ast}(y+z)-v^{\ast} (y)\qquad\qquad\qquad
$$
$$
	\qquad\quad\qquad
	-\I \left\langle z,\nabla_y v^{\ast} (y) \right\rangle ] dq(z)-a(y)I-g(y)\leq \mu \qquad y\in {\bf R^N}, 
$$
$$
	d+\sup_{\a\in \cal{A}} \{ \la \beta(\a),  \nabla_y v_{\ast} (y) \ra-c(\a)\}
	- a(y) \int_{{\bf R^N}} [v_{\ast} (y+z)-v_{\ast} (y)\qquad\qquad\qquad
$$
$$
	\qquad\quad\qquad
	-\I \left\langle z,\nabla_y v_{\ast} (y) \right\rangle ] dq(z)-a(y)I-g(y)\geq -\mu \qquad y\in {\bf R^N}. 
$$
The uniqueness of the nomber $d$ can be shown by the standard method (see \cite{al}). 
Thus, the claim was proved. \\
(ii) The existence of the number $d$ in (\ref{ae2}) can be shown in the same way to (i), by using Theorem 4.5.\\

{\bf Proof of Proposition 5.1. }
 (i) Put $\beta(\a)$$=b(x,\a)$ in (\ref{ae2}). Then, from Lemma 5.2, the statement follows. \\
(ii) In (\ref{ae1}), put $a(y)\equiv a$ ($a>0$ is a constant), $\beta(\a)=b(x,\a)$, $c(\a)=\left\langle b(x,\a),p \right\rangle $. The claim follows from Lemma 5.2. \\
(iii) In (\ref{ae1}), put $\beta(\a)=0$, $c(\a)=0$. From Lemma 5.2, the claim follows. \\

\section{Homogenizations.} 

$\quad$First, we confirm that the effective integro-differential operator has the uniform subellipticity. \\

(Uniform subelliptic operator) An integro-differential operator ${\mathcal I}(x,p,I)$ defined in ${\bf \Omega}\times$${\bf R^N}\times$${\bf R}$ is uniformly subelliptic if there exists $\theta>0$ such that 
\begin{equation}\label{subelliptic}
	{\mathcal I}(x,p,I+I')  \leq {\mathcal I}(x,p,I) - \theta I' \qquad \forall I'>0,\quad \forall (x,p,I)\in {\bf \Omega}\times{\bf R^N}\times{\bf R}. 
\end{equation}

{\bf Proposition 6.1.$\quad$} 
\begin{theorem}  
(i) Let $d_{x,p,I}$ be given by  (\ref{cell1}) in the case I, by (\ref{cell2}) 
 in the case II, and by (\ref{cell3}) in the case III, in the  sense of 
 (\ref{pw}). Put $\oI(x,p,I)=-d_{x,p,I}$ for any 
 $(x,p,I)\in \Omega \times {\bf R^N}\times {\bf R}$. Then,  $\oI$ is continuous and uniformly subelliptic. \\

(ii)  Let $d_{x,p,I}$ be given by  (\ref{qcell1}) in the case I', by (\ref{qcell2}) 
 in the case II', and by (\ref{qcell3}) in the case III', in the sense of 
 (\ref{pw}). Put $\oI(x,p,I)=-d_{x,p,I}$ for any 
 $(x,p,I)\in \Omega  \times {\bf R^N}\times {\bf R}$. Then,  $\oI$ is continuous and uniformly subelliptic. \\
\end{theorem}

{\bf Proof.}  (i) We use the perturbed test function method in \cite{ev1}. 
 We  prove the claim for the case II. The other cases can be treated similarly.  
 Take an arbitrary positive number $I'>0$. From (\ref{cell2}) and (\ref{pw}), for $\rho>0$  there exist bounded functions $v^{I}$ and $v^{I+I'}$ which satisfy
$$
	d(x,p,I)+ \sup_{\a\in \cal{A}} \{
 	\left\langle -b(x,\a), p+ \nabla v^{I}(y) \right\rangle \}  - 
	a(y) \int_{{\bf R^N}} [ v^{I}(y+z)-
	 v^{I}(y)
$$
$$
	- \I \left\langle z, \nabla  v^{I}(y) \right\rangle] dq(z) 
	-g_M(\gamma_1^{-1}y,...,\gamma_M^{-1}y)-a(y) I \leq  \rho \qquad y \in {\bf R^N}, 
$$
$$
	d(x,p,I+I')+ \sup_{\a\in \cal{A}} \{
 	\left\langle -b(x,\a), p+ \nabla v^{I+I'}(y) \right\rangle \}  - 
	a(y) \int_{{\bf R^N}} [ v^{I+I'}(y+z)- v^{I+I'}(y)
$$
\begin{equation}\label{super}
	-\I \left\langle z, \nabla  v^{I+I'}(y) \right\rangle] dq(z)
	-g_M(\gamma_1^{-1}y,...,\gamma_M^{-1}y) -a(y) (I+I') \geq - \rho \quad y \in {\bf R^N}. 
\end{equation}

By adding a constant if necessary, we may assume that $v^{I+I'}<v^{I}$. 
Let us prove (\ref{subelliptic}) for $\theta=a_0$, where $a_0>0$ is given in (\ref{a}) .   Assume that for fixed $x$, $p$, $I$ and $I'$,  there exists a constant $l>0$ such that 
\begin{equation}\label{contradiction}
	\oI(x,p,I+I') > \oI(x,p,I)- a_0 I' +l, 
\end{equation}
and we  look for a contradiction. 
We claim that $v^{I+I'}$ is a viscosity supersolution of 
$$
	-\oI(x,p,I)+ \sup_{\a\in \cal{A}} \{
 	\left\langle -b(x,\a), p+ \nabla v^{I+I'}(y) \right\rangle \} 
	 - a(y) \int_{{\bf R^N}} [ v^{I+I'}(y+z)-
	 v^{I+I'}(y)
$$
$$
	- \I \left\langle z, \nabla  v^{I+I'}(y) \right\rangle] dq(z)
	-g_M(\gamma_1^{-1}y,...,\gamma_M^{-1}y) -a(y) I \geq l - \rho \qquad y \in {\bf R^N}. 
$$
In fact, if there exists $\phi\in C^2$ such that $v^{I+I'}-\phi$ takes a global minimum at $y_0$, then from (\ref{super})
$$
	-\oI(x,p,I+I')+ \sup_{\a\in \cal{A}} \{
 	\left\langle -b(x,\a), p+ \nabla \phi (y_0) \right\rangle \} 
$$
$$
	 - \frac{a(y_0)}{2} \int_{|z|\leq \nu}  \left\langle (\nabla^2 \phi(y_0) + 2\delta I)z,z \right\rangle dq(z) 
	- a(y_0) \int_{|z| > \nu}
	[ \phi(y_0+z)-
	 \phi(y_0)
$$
$$
	- \I \left\langle z, \nabla  \phi(y_0) \right\rangle]  dq(z) 
	- g_M(\gamma_1^{-1}y_0,...,\gamma_M^{-1}y_0)  -a(y_0) (I+I') \geq  - \rho. 
$$
From (\ref{contradiction}),  
$$
	-\oI(x,p,I)+ \sup_{\a\in \cal{A}} \{
 	\left\langle -b(x,\a), p+ \nabla \phi (y_0) \right\rangle \} 
$$
$$
	 - 
	\frac{a(y_0)}{2}  \int_{|z|\leq \nu}  \left\langle (\nabla^2 \phi(y_0) + 2\delta I)z,z \right\rangle dq(z) 
	- a(y_0) \int_{|z| > \nu}
	[ \phi(y_0+z)-
	 \phi(y_0)
$$
$$
	- \I \left\langle z, \nabla  \phi(y_0) \right\rangle] dq(z)
	- g_M(\gamma_1^{-1}y_0,...,\gamma_M^{-1}y_0)  -a(y_0) I \geq  l- \rho, 
$$
which shows that $v^{I+I'}$ is the supersolution of the problem. \\
Since $v^{I}$, $v^{I+I'}$ are bounded, for $\l>0$ small enough, the following hold. 
$$
	\l v^{I} + \sup_{\a\in \cal{A}} \{
 	\left\langle -b(x,\a), p+ \nabla v^{I}(y) \right\rangle \}  
	- a(y) \int 
	[ v^{I}(y+z)-
	 v^{I}(y) \qquad\qquad\qquad
$$
$$
	- \I \left\langle z, \nabla  v^{I}(y) \right\rangle] dq(z)
	- g_M(\gamma_1^{-1}y,...,\gamma_M^{-1}y) -a(y) I \leq  2 \rho \qquad y \in {\bf R^N}, 
$$
$$
	\l v^{I+I'} + \sup_{\a\in \cal{A}} \{
 	\left\langle -b(x,\a), p+ \nabla v^{I+I'}(y) \right\rangle \}  
	- a(y) \int 
	[ v^{I+I'}(y+z)-
	 v^{I+I'}(y)\qquad
$$
$$
	- \I  \left\langle z, \nabla  v^{I+I'}(y) \right\rangle] dq(z) 
	-g_M(\gamma_1^{-1}y,...,\gamma_M^{-1}y)   -a(y) I \geq  l- 2 \rho \qquad y \in {\bf R^N}. 
$$
From the comparison principle, the above leads $\l(v^{I} - v^{I+I'} )\leq -l+4\rho$. Thus, by taking $\rho=\frac{l}{4}$ we get a contradiction to the fact that 
$v^{I+I'}<v^{I}$. Therefore, (\ref{subelliptic}) holds for $\theta=a_0$. \\
(ii) As shown in \S 5, for each case of I', II', and III', the integro-differential operator $\oI(x,p,I)$ is the limit of the sequence of operators 
$\oI_M(x,p,I)$ ($M=1,2,...$).  From (i), 
$$
	\oI_M(x,p,I+I') \leq  \oI_M(x,p,I) - a_0 I' \quad \forall I'>0, \quad \forall (x,p,I)\in {\bf R^N\times R^N \times R}. 
$$
 Therefore, (\ref{subelliptic})  holds for $\oI(x,p,I)$$=\lim_{M\to \infty} \oI_M(x,p,I)$, too. \\

$\qquad$The following is the main result of the quasi-periodic homogenization. \\

{\bf Theorem 6.2.$\quad$} 
\begin{theorem} Assume that (\ref{b}), (\ref{a}), (\ref{nonreso}), (\ref{holder}) and  (\ref{gholder}) hold.
 If $\a=1$ and $\beta(x,\a)\not \equiv 0$, assume that $a(y)\equiv a$ ($a>0$ is a constant). 
Let $\ue$ be the solution of (\ref{quasi})-(\ref{bc}). Then, there exists a function $\ou$ such that
$$
	\lim_{\varepsilon\to 0} \ue (x)= \ou(x) \qquad \hbox{uniformly in}\quad x\in {\bf R^N}. 
$$
The function $\ou$ is the unique bounded solution of (\ref{effect})-(\ref{bc}) with the effective integro-differential operator $\oI$, given by 
$$
	\oI(x,p,I)=-d_{x,p,I}\quad  \hbox{for any}  \quad (x,p,I)\in {\bf \Omega \times R^N \times R}.
$$ 
 The right-hand side $d_{x,p,I}$ is given by Proposition 4.6. \\
\end{theorem}

{\bf Proof.} It is enough to prove the case of $\a=1$. The other cases can be shown similarly. \\
 (Step 1) First, remark that from Proposition 4.6, the effective integro-differential equation (\ref{effect}) is well-definded. 
From Proposition 6.1, $\oI$ is subelliptic, and the comparison principle holds for  (\ref{effect})-(\ref{bc}) (see \cite{ar3},\cite{ar5}, \cite{ar8}, \cite{bbp} and \cite{bi} for example).  \\
 (Step 2) From the comparison for (\ref{quasi})-(\ref{bc}), there exists a constant $M>0$ such that $|\ue|\leq M$ for any $\varepsilon\in (0,1)$. Therefore, we can take 
$$
	u^{\ast}(x)=\overline{\lim}_{\varepsilon \to 0,x'\to x} \ue (x'), \quad u_{\ast}(x)=\overline{\lim}_{\varepsilon \to 0,x'\to x} \ue (x'). 
$$
We prove that $u^{\ast}$ is a viscosity subsolution of (\ref{effect}) by the argument by the contradiction. 
So, assume that $u^{\ast}$ is not a subsolution of (\ref{effect}): there exists a function $\phi\in C^2$ such that $u^{\ast}-\phi$ takes a global strict maximum at a  point $\hx$$\in {\bf R^N}$, $u^{\ast}(\hx)= \phi(\hx)$, and 
$$
	 \phi(\hx)  + 
 	\oI(\hx,\nabla \phi(\hx), I[\phi](\hx))= 3\gamma >0, 
$$
where $\gamma>0$ is a constant. 
 From the continuity of $\oI$, for $r>0$ small enough 
$$
	 \phi(x) + 
 	\oI(x,\nabla \phi(x), I[\phi](x)) \geq  \gamma  \quad \hbox{in}\quad B_r(\hx), 
$$
where $B_r(\hx)$$=\{x|\quad |x-\hx|<r \}$. 
From the definition of  $\oI(\hx, \nabla \phi(\hx), I[\phi](\hx)) $, for the above $\gamma >0$ there exist $\uv$, $\ov$ such that 
$$
	-\oI(\hx,\nabla \phi(\hx),I[\phi](\hx))+ \sup_{\a\in \cal{A}} \{
 	\left\langle -b(\hx,\a),  \nabla \uv(y) \right\rangle \}  - a(y) \int [ \uv (y+z)- \uv (y)
$$
\begin{equation}\label{uv}
	- \I \left\langle z, \nabla  \uv (y) \right\rangle] \frac{1}{|z|^{N+\a}}dz
	-g_M(\gamma_1^{-1}y,...,\gamma_M^{-1}y)  -a(y) I [\phi](\hx) \leq  \gamma \quad y \in {\bf R^N}, 
\end{equation}
$$
	-\oI(\hx,\nabla \phi(\hx),I[\phi](\hx)+ \sup_{\a\in \cal{A}} \{
 	\left\langle -b(\hx,\a),  \nabla \ov(y) \right\rangle \}  - a(y) \int [ \ov (y+z)- \ov (y)
$$
\begin{equation}\label{ov}
	- \I \left\langle z, \nabla  \ov (y) \right\rangle] \frac{1}{|z|^{N+\a}}dz
	-g_M(\gamma_1^{-1}y,...,\gamma_M^{-1}y)   -a(y) I[\phi](\hx) \geq  - \gamma \quad y \in {\bf R^N}. 
\end{equation}
We can assume that $\uv$, $\ov$ are Lipschitz continuous, for if not we regularize them by the sup and the inf convolutions respectively. 
Put $\phi_{\varepsilon}=\phi(x)+\varepsilon_1^{\a}\ov(\frac{x}{\varepsilon_1})$. We claim that $\phi_{\varepsilon}$ is a viscosity supersolution of 
$$
	\phi_{\varepsilon}+ \sup_{\a\in \cal{A}} \{
 	\left\langle -b(x,\a), \nabla  \phi_{\varepsilon} \right\rangle \} 
	- a(\frac{x}{\varepsilon_1}) \int_{{\bf R^N}} [\phi_{\varepsilon}(x+z)-\phi_{\varepsilon}(x)\qquad\qquad\qquad
$$
\begin{equation}\label{claim}
	\qquad
	-\I \left\langle z, \nabla \phi_{\varepsilon}(x) \right\rangle] \frac{1}{|z|^{N+\a}}dz
	-g_M(\frac{x}{\varepsilon_1},...,\frac{x}{\varepsilon_M}) \geq \gamma \qquad \hbox{in} \quad B_r(\hx), 
\end{equation}
for $r>0$ small enough.  To see this, assume that for $\psi\in C^2$ $\phi_{\varepsilon}-\psi$ attains its minimum at $\ox\in U_r(\hx)$, and that 
 $\phi_{\varepsilon}(\ox)=\psi(\ox)$. We are to show 
$$
	\phi_{\varepsilon}(\ox)+ \sup_{\a\in \cal{A}} \{
 	\left\langle -b(\ox,\a), \nabla  \psi(\ox) \right\rangle \} 
	- a(\frac{\ox}{\varepsilon_1})
	\int_{{\bf R^N}} [\psi (\ox+z)-\psi(\ox)- 
$$
$$
	\I \left\langle z, \nabla \psi(\ox) \right\rangle] \frac{1}{|z|^{N+\a}}dz
	-g_M(\frac{\ox}{\varepsilon_1},...,\frac{\ox}{\varepsilon_M})  \geq \gamma \qquad \hbox{in} \quad B_r(\hx). 
$$
Put $\beta(y)=\frac{1}{\varepsilon_1^{\a}}(\psi-\phi)(\varepsilon_1 y)$. Since $(\ov-\beta)(y)$ attains its minimum at $\oy=\frac{\ox}{\varepsilon_1}$, and since  
$\ov$ is the viscosity supersolution of (\ref{ov}), 
$$
	-\oI(\hx,\nabla \phi(\hx),I[\phi](\hx)) 
	+\sup_{\a \in {\mathcal A}} \{ 
	\left\langle -b(\hx,\a), \nabla  (\psi-\phi)(\ox) \right\rangle 
	\}\qquad\qquad\qquad\qquad
$$
$$- a(\frac{\ox}{\varepsilon_1})
	\int_{{\bf R^N}} [ \frac{\psi- \phi }{\varepsilon_1^{\a}}(\ox+\varepsilon_1 z)
	-\frac{\psi- \phi }{\varepsilon_1^{\a}}(\ox)
	- \I \left\langle \varepsilon_1 z,  \frac{\nabla (\psi- \phi) }{\varepsilon_1^{\a}}(\ox) \right\rangle] \frac{1}{|z|^{N+\a}}dz
$$
$$\qquad\qquad\qquad\qquad\qquad\qquad\qquad
	- g_M(\frac{\ox}{\varepsilon_1},...,\frac{\ox}{\varepsilon_M}) -a(\frac{\ox}{\varepsilon_1}) I[\phi](\hx) \geq  -\gamma. 
$$
From the continuity of $\oI$ (Proposition 6.1 (i)), 
$$
	-\oI(\ox,\nabla \phi(\ox),I[\phi](\ox)) - a(\frac{\ox}{\varepsilon_1})
	\int_{{\bf R^N}} [    (\psi- \phi ) (\ox+z)-(\psi- \phi )(\ox)\qquad\qquad\qquad\qquad
$$
$$\qquad
	- \I \left\langle \varepsilon_1 z, \nabla  (\psi- \phi)(\ox) \right\rangle] \frac{1}{|z|^{N+\a}}dz
	-g_M(\frac{\ox}{\varepsilon_1},...,\frac{\ox}{\varepsilon_M}) -a(\frac{\ox}{\varepsilon_1}) I[\phi](\hx) \geq - \gamma. 
$$
Therefore,  for $r>0$ and $\varepsilon>0$ small enough 
$$
	\phi_{\varepsilon}(\ox)+ \sup_{\a\in \cal{A}} \{
 	\left\langle -b(\ox,\a), \nabla  \psi(\ox) \right\rangle \} 
	- a(\frac{\ox}{\varepsilon_1})
	\int_{{\bf R^N}} [\psi (\ox+z)-\psi (\ox) \qquad\qquad\qquad\qquad
$$$$\qquad\qquad\qquad\qquad\qquad
	 -\I \left\langle z, \nabla \psi (\ox) \right\rangle] \frac{1}{|z|^{N+\a}}dz
	-g_M(\frac{\ox}{\varepsilon_1},...,\frac{\ox}{\varepsilon_M})
$$
$$
	\geq 
	\phi_{\varepsilon}(\ox)+ \sup_{\a\in \cal{A}} \{
 	\left\langle -b(\ox,\a), \nabla  \phi (\ox) \right\rangle \} 
	+ \gamma+ \oI(\ox,\nabla \phi(\ox),I[\phi](\ox))
	- a(\frac{\ox}{\varepsilon_1})
	\int_{{\bf R^N}} [\phi (\ox+z)-
$$$$\qquad\qquad\qquad\qquad
	\phi (\ox) - \I \left\langle z, \nabla \phi (\ox) \right\rangle] \frac{1}{|z|^{N+\a}}dz
	+ a(\frac{\ox}{\varepsilon_1})I[\phi](\hx)\geq 2\gamma, 
$$
where we used $\phi_{\varepsilon}(\ox)=\psi(\ox)$, and the fact that $v$ is Lipschitz continuous.  
Hence, (\ref{claim}) was proved.  From  the comparison $u_{\varepsilon} \leq \phi_{\varepsilon}-2\gamma$,  and since $\gamma>0$ is arbitrary $u_{\varepsilon} \leq \phi_{\varepsilon}$. 
Therefore, 
$$
	u^{\ast}(x)= \overline{\lim_{\varepsilon\to 0,x'\to x}} u_{\varepsilon}(x')\leq \lim_{\varepsilon \to 0} \phi_{\varepsilon}(x) \quad  \hbox{in} \quad B_{r}(\hx). 
$$ 
However, this contradicts to the fact that $u^{\ast}-\phi$ takes its strict maximum at $\hx$, and we have proved that $u^{\ast}$ is a viscosity subsolution of (\ref{effect}). 
In parallel, we can prove that  $u_{\ast}$ is a viscosity supersolution of (\ref{effect}).\\
(Step 3) As we have confirmed in  Step 1, the comparison principle holds for (\ref{effect})-(\ref{bc}). The fact that $u^{\ast}$ and $u_{\ast}$ are respectively a subsolution and a supersolution  
of  (\ref{effect}) leads 
$u^{\ast}\leq u_{\ast}$.  At the same time, from the definition $u_{\ast}\leq u^{\ast}$. Therefore, $u^{\ast}=u_{\ast}$, and $\lim_{\varepsilon\to 0}\ue=\ou$ exists. From Step 2, $\ou$ is the unique solution of (\ref{effect})-(\ref{bc}).\\

The following is our main  result of the almost periodic homogenizations. \\

{\bf Therem 6.3.$\quad$} 
\begin{theorem} Assume that (\ref{b}), (\ref{a}), (\ref{holder}) and (\ref{gholder}) hold, and that $g$ is uniformly almost periodic in the sense of Bohr in ${\bf R^N}$.  
 If $\a=1$ and $\beta(x,\a)\not \equiv 0$, assume that $a(y)\equiv a$ ($a>0$ is a constant). 
Let $\ue$ be the solution of (\ref{almostex})-(\ref{bc}). Then, there exists a function $\ou$ such that
$$
	\lim_{\varepsilon\to 0} \ue (x)= \ou(x) \qquad \hbox{uniformly in}\quad x\in {\bf R^N}. 
$$
The function $\ou$ is the unique  bounded solution of (\ref{effect})-(\ref{bc}) with the effective integro-differential operator $\oI$, given by 
$$
	\oI(x,p,I)=-d_{x,p,I}\quad  \hbox{for any}  \quad (x,p,I)\in {\bf \Omega \times R^N \times R}.
$$ 
 The right-hand side $d_{x,p,I}$ is given by Proposition 5.1. \\
\end{theorem}

{\bf Proof.} From Lemma B, we can take a sequence of  functions $g_M(y_1,...,y_M)$ ($M=1,2,...$) periodic in ${\bf T^{MN}}$,  and a sequence of positive numbers $\{\gamma_i\}$ ($i=1,2,...$) satisfying the non-resonance condition (\ref{nonreso}), such that 
\begin{equation}\label{conv}
	g(y)= \lim_{M\to \infty} g_M(\gamma_1^{-1} y,...,\gamma_M^{-1} y)    \quad \hbox{uniformly in }\quad y\in {\bf R^N}. 
\end{equation}
We assume that $\gamma_1=1$. Let $\varepsilon_1>0$, and put $\varepsilon_i=\gamma_i \varepsilon_1$ for any $2\leq i\leq M$. 
For $\varepsilon=(\varepsilon_1,...,\varepsilon_M)$, let $u_{\varepsilon}^{M}$ be the solution of 
\begin{equation}\label{Mquasi}
	\ue^M + \sup_{\a\in \cal{A}} \{
 	\left\langle -b(x,\a), \nabla\ue^M \right\rangle \} - a(\frac{x}{\varepsilon_1})\int_{{\bf R^N}} [\ue^M(x+z)-\ue^M(x)
\end{equation}
$$
	- \I \left\langle z, \nabla\ue^M(x) \right\rangle] \frac{1}{|z|^{N+\a}}dz
	-g_M(\gamma_1^{-1} \frac{x}{\varepsilon_1},...,\gamma_M^{-1} \frac{x}{\varepsilon_1})=0 \qquad x\in {\bf R^N}, 
$$
and (\ref{bc}). 
Then, by comparing the above equation with (\ref{almostex}), there exists a sequence of constants $c_M>0$ such that $\lim_{M\to \infty} c_M =0$, 
\begin{equation}\label{eM}
	u_{\varepsilon}^{M}-c_M \leq \ue \leq u_{\varepsilon}^{M}+c_M \quad \forall \varepsilon>0, \quad \forall M\in {\bf N}. 
\end{equation}
	From Theorem 6.2, there exists a sequence of integro-differential operators $\oI_M$ such that the limit $\lim_{\varepsilon \to 0}\ue^M(x)$$=u^{M}(x)$ 
satisfies 
$$
	u^{M}(x)+ \oI_M (x,\nabla u^{M}(x), I[u^{M}](x))=0 \quad \hbox{in}\quad {\bf R^N}, 
$$
and (\ref{bc}). 
Put ${u}^{\ast}(x)=\lim_{M\to \infty,y\to x} u^{M}(y)$,  ${u}_{\ast}(x)=\lim_{M\to \infty,y\to x} u^{M}(y)$. Then, ${u}^{\ast}$ and $u_{\ast}$ respectively satisfy 
$$
	{u}^{\ast} + \oI (x,\nabla {u}^{\ast}, I[{u}^{\ast}](x)\leq 0 \quad \hbox{in}\quad {\bf R^N}, 
$$
$$
	{u}_{\ast} + \oI (x,\nabla {u}_{\ast}, I[{u}_{\ast}](x)\geq 0 \quad \hbox{in}\quad {\bf R^N}, 
$$
where $\oI(x,p,I)=\lim_{M\to \infty} \oI_M (x,p,I)$. From the comparison, we get ${u}_{\ast} \leq {u}^{\ast} \leq {u}_{\ast}$. Thus, 
$$
	\lim_{M\to \infty} u^M = \exists \ou(x) \quad \hbox{in}\quad {\bf R^N}. 
$$
Now, we first let $\varepsilon\to 0$ in (\ref{eM}) to have
$$
	u^{M}-c_M \leq \underline{\lim}_{\varepsilon\to 0}\ue \leq  \overline{\lim}_{\varepsilon\to 0}\ue
	\leq u^{M}+c_M \quad \forall M\in {\bf N}, 
$$ 
then let $M\to \infty$ to have 
$$
	\ou(x) \leq \underline{\lim}_{\varepsilon\to 0}\ue \leq  \overline{\lim}_{\varepsilon\to 0}\ue
	\leq \ou(x). 
$$ 
Thus, $\lim_{\varepsilon\to 0} u_{\varepsilon}=\ou$ exists which is the unique solution of (\ref{effect})-(\ref{bc}). \\



\end{document}